\documentclass[12pt]{article}
\usepackage{graphicx} 
\usepackage{epstopdf}
\usepackage{subfigure}
\usepackage{color}
\usepackage[english]{babel}
\usepackage{amsmath,amsthm}
\usepackage{amssymb}
\usepackage{amsfonts}
\usepackage{booktabs}
\usepackage{cite}

\usepackage{graphicx}
\usepackage{subfigure}
\usepackage{tabularx}
\title{On the classification and dispersability of circulant graphs with two jump lengths}

\author {  Xiaoxiang Yu, Zeling Shao, Zhiguo Li{$^*$}\\
{\small School of Science, Hebei University of Technology, Tianjin 300401, China}
\date{}
\footnote{Corresponding author. E-mail: zhiguolee@hebut.edu.cn}
}

\begin{document}
\baselineskip 0.65cm

\maketitle

\begin{abstract}

 In this paper, we give the classification of circulant graphs $C(\mathbb{Z}_{n},S)$ with $|S|=2$ and completely solve the dispersability of circulant graphs $C(\mathbb{Z}_{n},\{1, k\})$.


\noindent\textbf{2020 MR Subject Classification.} 05C10
\end{abstract}

\section{Introduction}
The research of book embedding is to find an optimal embedding that meets certain conditions. It can greatly optimize some performance in the field of computer science, including fault-tolerant computing, multilayer VLSI layout and so on(see $[1,2,3]$).
Let $\psi$ be a permutation for all vertices of a graph $G$.~A \emph{layout} $\Psi=(G,\psi)$ for $G$ is to arrange all vertices along a circle in the order $\psi$ and join the edges of $G$ as chords.~Let $S$ be a color set and $|S|=n$.~A triple $(G,\psi,c)$ is an $n$-page \emph{book embedding} if $c:E(G) \rightarrow S$ is an edge-coloring such that $c(e^{\prime})\neq c(e^{\prime\prime})$ when $e^{\prime}$ and $e^{\prime\prime}$ cross in $\Psi$.~The \emph{book thickness}~$bt(G)$ of $G$ is the minimum integer $n$ such that an $n$-page book embedding exists.~A book embedding $(G,\psi,c)$ is \emph{matching} if the edge-coloring $c$ is proper.~The \emph{matching book thickness} $mbt(G)$ of $G$ is the minimum integer $n$ such that an $n$-page matching book embedding exists.~We call $G$ \emph{dispersable} if $mbt(G)=\Delta(G)$ and \emph{nearly dispersable} if  $mbt(G)=\Delta(G)+1$(see $[4]$).

The dispersability of some families of graphs has been studied.~Overbay ${[5]}$ discussed the dispersability of complete bipartite graphs, even cycles, binary $n$-cubes$(n\geq1)$, trees, complete~graphs. The dispersability of the generalized Petersen graph and the pseudo-Halin graph is obtained in $[6],[7]$ respectively. Meanwhile,~the dispersability~of~the~Cartesian~product~of~two~cycles~has~some developments. Kainen ${[8]}$ showed that $mbt(C_{p}\Box C_{q})$ is $4$, when $p,q$ are both even and $mbt(C_{p}\Box C_{q})$ is $5$, when $p$ is even, $q$ is odd. Shao,~Liu,~Li~$[9]$~obtained~that~$mbt(K_{n}\Box C_{q})=\Delta(K_{n}\Box C_{q})+1$~for~$n,q \geq3$,~which~implies~$mbt(C_{3}\Box C_{q})=5$.~In~$[4]$~it was proved that $mbt(C_{5}\Box C_{n})=5$ for $n\geq 3$.~Recently,~the authors $[10]$ proved that $mbt(C_{2m+1}\Box C_{2n+1})=5(m,n \geq 3)$,~which completely solved the dispersability of the Cartesian product of two cycles.

Let $S$ be a subset of the group $\mathbb{Z}_{n}$ such that each element $k$ in $S$ satisfies $1\leq k\leq \lfloor n/2 \rfloor$. A \emph{circulant graph} $G=C(\mathbb{Z}_{n},S)$ is a graph whose vertex set is $ \mathbb{Z}_{n}$,~two vertices $v_i,v_j$ are adjacent if $v_i-v_j\equiv k\pmod{n}$ and $k\in S$. For convenience, unless otherwise expressly indicated, the group $\mathbb{Z}_n$ is the set $\{1,2,\cdots,n \}$ throughout this paper.~And each element of $S$ is called a jump length.~If $|S|=1$, it is easy to see $G$ is dispersable if $G$ is bipartite; otherwise, $G$ is nearly dispersable [11].

\vspace{-0.6em}
\begin{figure}[htbp]
\centering
\includegraphics[height=3.4cm, width=0.62\textwidth]{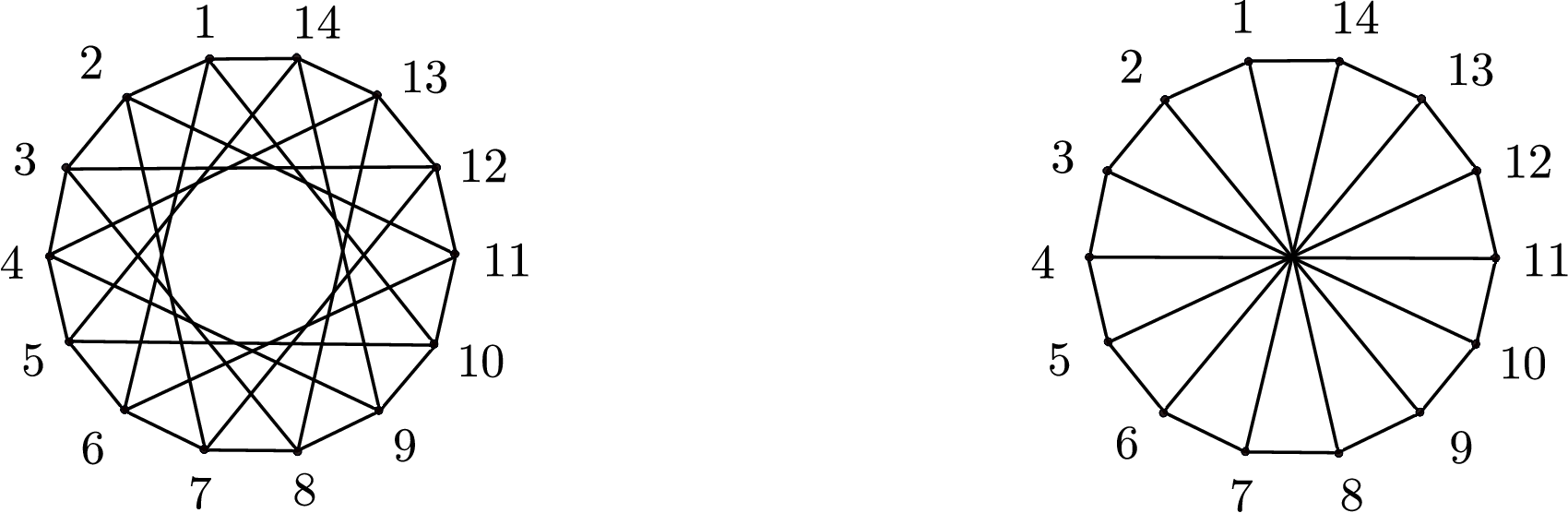}
\centerline{Fig.1  ~The circulant graph $C(14,5)$(left) and $C(14,7)$(right).}%
\end{figure}

\vspace{-1.2em}

Joslin, Kainen and Overbay $[11]$ gave some results of the dispersability of circulant graphs $C(\mathbb{Z}_{n},S)$, where $S$ is a subset of $\{1,2,3\}$ with size $2$. In addition, they showed if $n$ is a multiple of $2k+1$, then $C(\mathbb{Z}_{n},\{1,2,\cdots,k\})$ is nearly dispersable; if $n$ is a multiple of $12$, then $C(\mathbb{Z}_{n},\{1,2,3\})$ is nearly dispersable; if $2k\mid n~(k\geq 3)$, then $C(\mathbb{Z}_{n},\{1,3,5,\cdots,l\})$ is dispersable, where $l$ is the largest odd integer not exceeding $k$.

In this paper, We focus on the case of the set $S$ being any subset of $\mathbb{Z}_{n}$ with size $2$.~For convenience,~we denote $C(\mathbb{Z}_{n},\{1,k\})$ by $C(n,k)$ (see Fig.1 for $C(14,5)$(left) and $C(14,7)$(right)).~Obviously,~$C(n,k)\cong C(n,n-k)$ if $\lceil n/2 \rceil \leq k \leq n-1.$ We mainly give the classification of circulant graphs $C(\mathbb{Z}_{n},\{k_1,k_2\})$ and prove $C(n, k)$ is dispersable when $n$ is even and $k$ is odd; Otherwise, $C(n,k)$ is nearly dispersable.



The paper is organized as follows. We introduce some definitions and properties in Section 2. Section 3 gives the classification of circulant graphs $C(\mathbb{Z}_{n},\{ k_1, k_2\})$ with two jump lengths. In section 4, we obtain the dispersability of $C(n,k)$, $n$ is even.~For $n$ is odd, section 5 and section 6 discuss the cases of $gcd(n,k) \neq 1$ and $gcd(n,k)=1$ respectively.

\section{Preliminaries}
For convenience, we give some notations. Let $gcd(n,k)$ be the greatest common denominator of $n$ and $k$. If $X$ is an ordered vertex set, let $X^{-}$ be the reverse of $X$ and $|X|$ be the size of $X$. The vertices of $C(n,k)$ are in the sense of modulo $n$. In order to give the classification of circulant graphs $C(\mathbb{Z}_{n},\{k_1,k_2\})$, let us  recall the definition of Cartesian graph bundles.\\
\noindent
\textbf{Definition 2.1.~}$^{[12]}$~\emph{Let $B,F$ be graphs. A graph $G$ is a Cartesian graph bundle with fiber $F$ over the base graph $B$ if there is a graph map $p:G\rightarrow B$ such that for each vertex $v\in V(B),p^{-1}(v)\cong F,$ and for each edge $e\in E(B),p^{-1}(e)\cong K_2 \Box F.$ Let $\varphi:E(B)\rightarrow Aut(F)$ be a mapping which assigns an automorphism of the graph $F$ to any edge of $B$. The bundle $G$ is denoted by $G=B\Box^{\varphi}F$.}\\
\indent
Automorphisms of a cycle $C_t$ are of two types~$[13]$. A cyclic shift of the cycle by $d(0\leq d\leq t)$ elements is called the cyclic $d$-shift and other automorphisms of $C_t$ are  called reflections(see Fig.2).

\vspace{-0.6em}
\begin{figure}[htbp]
\centering
\includegraphics[height=2.6cm, width=0.66\textwidth]{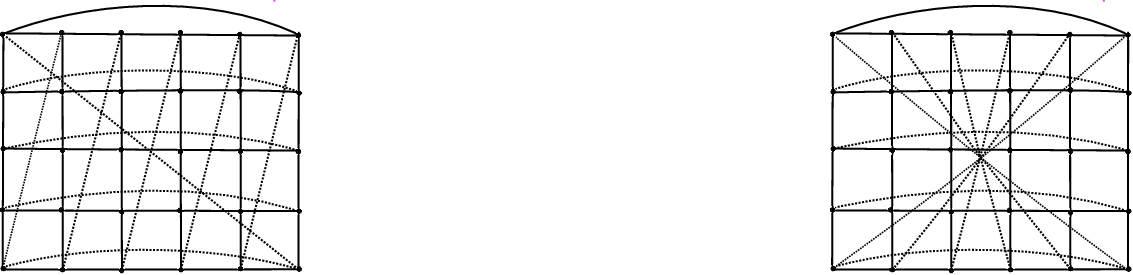}
\centerline{Fig.2  ~The Cartesian graph bundles $C_5\Box^{\varphi}C_6$, where $\varphi$ is a 1-shift(left) and a reflection(right).}%
\end{figure}
\vspace{-1.2em}

Next, we recall some properties of matching book embedding and  circulant graphs $C(n,k)$.

\vspace{0.0em}
\noindent \textbf{Lemma 2.1.~}$^{[5]}$\emph{~If a regular graph $G$ is dispersable, then $G$ is bipartite.}

\vspace{0.0em}
\noindent \textbf{Lemma 2.2.~}$^{[5]}$\emph{
For any simple graph $G$, we have $\Delta(G)\leq mbt(G)$ .}

\vspace{0.0em}
\noindent \textbf{Lemma 2.3.~}$^{[14]}~$\emph{$C(n,k)$ is bipartite if and only if $n$ is even and $k$ is odd.}

\vspace{0.0em}
\noindent \textbf{Lemma 2.4.~}$^{[14]}~$\emph{The maximum degree of C(n,k) is as follows.}
\vspace{-0.7em}
$$\Delta(C(n,k))=\left\{\begin{array}{l}
3,\ \ \ \ \ \ \ \ \ n ~is ~even, ~k=n/2;\\
4,\ \ \ \ \ \ \ \ \ \ \ \ \ \ \ \ \ \ \ \ \ \ \ \ \ \ \ \ \ else.\\
\end{array} \right.$$
\vspace{-2em}

According to the above lemmas, the following result holds.

\noindent \textbf{Lemma 2.5.~}\emph{The matching book thickness of C(n,k) has the following result.}
\vspace{-0.6em}
$$mbt(C(n,k))\geq\left\{\begin{array}{l}
\Delta(C(n,k)),\ \ \ \ \ \ \ \ n ~is ~even, ~k~is~odd;\\
\Delta(C(n,k))+1,\ \ \ \ \ \ \ \ \ \ \ \ \ \ \ \ \ \ \ \ \ \ \ else.\\
\end{array} \right.$$
\vspace{-2em}

In order to calculate the matching book thickness of $C(n,k)$, the following definition, lemmas and notations are frequently used to give vertex orderings and edge colorings.\\
\noindent \textbf{Lemma 2.6.~}$^{[11]}~$\emph{$C(n,k)(k\neq n/2)$ is the edge disjoint union of $d+1$ cycles, where $d=gcd(n,k)$.}\\
\noindent \textbf{Proof.~}Assume~$C_0=(V,E)$~is~a~cycle,~where~$V=\mathbb{Z}_{n},E=\{(i,j)~|~(i-j)\equiv 1(\text{mod}~n)\}$.~Let $C_i(1 \leq i \leq d)$ be a cycle with the vertex set $V(C_i)=$ $\{i, i+k, i+2k, \cdots, i+(n/d-1)k\}(1 \leq i \leq d)$ such that two vertices $v_i,v_j$ are adjacent if $v_i-v_j \in \{\pm k\}(\text{mod}~n)$.~It is easy to see that $C_i$ is a subgraph of $C(n,k)$ for $i=0,1,\cdots,d$, and all edges of these $d+1$ cycles are disjoint and account for all the edges of $C(n,k)$.
\hfill{$\square$}

\noindent \textbf{Definition 2.2.~}\emph{Given a layout $\Psi=(G,\omega)$ of a graph $G$, let $S$ be a color set of size $mbt(G)$.}
\emph{We say the graph $G$ can be colored well in the order $\omega$ if there is a proper edge coloring $c:E(G)\rightarrow S$ such that $c(e)\neq c(e^{\prime})$ when $e$ and $e^{\prime}$ cross in $\Psi$ for $e,e^{\prime}\in E(G)$.}

\noindent \textbf{Remark 1.~}If $G$ is colored well,~it is easy to see that an $mbt(G)$-page matching book embedding of $G$ is obtained.

\vspace{0.0em}
\noindent \textbf{Lemma 2.7.~}\emph{Let $n,k(n>k\geq2)$ be positive integers, $n=ak+r$.~Assume $gcd(n,k)=1$, then there exists a finite nonnegative integer $m$ in the following algorithm so that $r_{m+1}=0, k_{m+1}=1$.}
\vspace{-1.7em}
\begin{equation}
n=ak+r
\end{equation}
\vspace{-1.7em}
\vspace{-0.6em}
\begin{equation}
k=a_1(k-r)+r_1=a_1k_1+r_1
\end{equation}
\vspace{-1.7em}
\vspace{-0.6em}
\begin{equation}
k_1=a_2(k_1-r_1)+r_2=a_2k_2+r_2
\end{equation}
\vspace{-1.7em}
\vspace{-0.6em}
\begin{equation}
k_2=a_3(k_2-r_2)+r_3=a_3k_3+r_3
\end{equation}
\vspace{-1.7em}
\vspace{-0.8em}
$$\vdots$$
\vspace{-2em}
\vspace{-0.8em}
\begin{equation}\tag {$m$+1}
k_{m-1}=a_m(k_{m-1}-r_{m-1})+r_m=a_mk_m+r_m
\end{equation}
\vspace{-1.5em}
\vspace{-0.8em}
\begin{equation}\tag {$m$+2}
k_m=a_{m+1}(k_m-r_m)+r_{m+1}=a_{m+1}k_{m+1}+r_{m+1}
\end{equation}

\vspace{-0.4em}
\noindent \textbf{Proof.} ~Without loss of generality, we denote $k,r$ by $k_0,r_0$ respectively.
Firstly, it is easy to see that  the sequence$\{r_i\}_{i=1}^{m+1}$ is monotonically decreasing. In fact, $k_{i-1}-r_i=a_ik_i\geq k_i=k_{i-1}-r_{i-1}$ for any $a_i\geq 1$, hence $r_i\leq r_{i-1}$. If there is an integer $j(1\leq j\leq m)$ such that $r_j=r_{j-1}>0$, we write $(j+1)$ as $k_{j-1}=1*(k_{j-1}-r_{j-1})+r_j$. If $r_j~|~k_{j-1}$, after at most $k_{j-1}/r_j$ steps, the remainder is zero; otherwise, there exists an integer $l(j<l\leq m)$ such that $r_l<r_j$. So there exists a minimum finite nonnegative integer $m$ such that $r_{m+1}=0.$

Next, before proving $k_{m+1}=1$,~we  prove that $gcd(k_{m},r_{m})=1$. Assume $gcd(k_{m},r_{m})=s\geq 2$.
According to the equality $(m+1)$, we have $gcd(k_{m-1},k_{m})=gcd(k_{m-1},k_{m-1}-r_{m-1})=gcd(k_{m-1},r_{m-1})=s.$~
In a similar way, we can conclude that $gcd(k_{i},r_{i})=gcd(k,r)=gcd(n,k)=s, 1 \leq i \leq m$, which is a contradiction. We simplify $(m+2)$ to the following:
\vspace{-0.8em}
$$\cfrac{(a_{m+1}-1)}{a_{m+1}}=\cfrac{r_m}{k_m}.$$

\vspace{-0.6em}
\noindent Since $gcd(k_{m},r_{m})=1$, $k_{m+1}=k_m-r_m=1$.
\hfill{$\square$}

\vspace{0.3em}
Assume $G=C(n,k)$ and $n=ak+r$. If $gcd(n,k)=1$, in order to give the vertex ordering of the optimal matching book embedding of $G$, we construct some ordered vertex sets as follows.\\
\indent
Let $P_i(1 \leq i \leq t)$ be ordered vertex sets, specifically,

\indent
$(i)$~If $r=1$, let $t=k$(see Example $(i)$) and
\vspace{-0.8em}
$$P_i=\left\{\begin{array}{l}
\{i,i+k,i+2k,\cdots,i+(a-1)k,i+ak\},\ \ \ \ \ \ \ \ \ \ \ \ \ \ \ \ \ \ i=1;\\
\{i,i+k,i+2k,\cdots,i+(a-1)k\},\ \ \ \ \ \ \ \ \ \ \ \ \ \ \ \ \ \ \ \ \ \ 2 \leq i \leq k.\\
\end{array} \right.$$
\vspace{-2em}

\indent
$(ii)$~If $r=k-1$, let $t=k$(see Example $(ii)$) and
\vspace{-0.8em}
$$P_i=\left\{\begin{array}{l}
\{i,i+k,i+2k,\cdots,i+(a-1)k,i+ak\},\ \ \ \ \ \ 1 \leq i\leq k-1;\\
\{i,i+k,i+2k,\cdots,i+(a-1)k\},\ \ \ \ \ \ \ \ \ \ \ \ \ \ \ \ \ \ \ \ \ \ \ \ \ \ i=k.\\
\end{array} \right.$$
\vspace{-2em}

$(iii)$~According to Lemma 2.7, there exists a positive integer $m$ such that $k_{m+1}=1,r_{m+1}=0$, if $1<r <k-1$. If there exists a minimum positive integer $s(s<m)$ such that $k_s \geq 3,r_s=1$ in the algorithm, let $t=k_{s}+1$(see Example $(iii)(iv)$) and, for convenience, we consider $C(n,n-k)$,
\vspace{-1.2em}
$$P_i=\left\{\begin{array}{l}
\{i,i+(n-k),i+2(n-k),\cdots,(i+1)-(n-k)\},\ \ \ \ 1 \leq i \leq k_s;\\
\{i,i+(n-k),i+2(n-k),\cdots,1-(n-k)\},\ \ \ \ \ \ \ \ \ \ \ \ i=k_s+1.\\
\end{array} \right.$$
\vspace{-1.6em}

\noindent
Otherwise, let $t=k_{m}$(see Example $(v)$) and, for convenience, we consider $C(n,k)$,
\vspace{-0.6em}
$$P_i=\left\{\begin{array}{l}
\{i,i+k,i+2k,\cdots,(i+1)-k \},\ \ \ \ \ \ \ \ \ \ \ \ \ \ \ \ \ 1 \leq i \leq k_{m}-1;\\
\{i,i+k,i+2k,\cdots,1-k \},\ \ \ \ \ \ \ \ \ \ \ \ \ \ \ \ \ \ \ \ \ \ \ \ \ \ \ \ \ \ \ \ \ \ \ i=k_{m}.\\
\end{array} \right.$$
\vspace{-2em}
\vspace{0.2em}

\noindent
\noindent \textbf{Example. }$(i)$~$C(56,5)$: $a=11, r=1.$\\
\noindent
$P_1$=$\{1,6,11,16,21,26,31,36,41,46,51,56\};$ $P_2$=$\{2,7,12,17,22,27,32,37,42,47,52\};$\\
\noindent
$P_3$=$\{3,8,13,18,23,28,33,38,43,48,53\};$ $P_4$=$\{4,9,14,19,24,29,34,39,44,49,54\};$\\
\noindent
$P_5$=$\{5,10,15,20,25,30,35,40,45,50,55\}.$\\
\noindent
$(ii)$~$C(53,9)$: $a=5, r=8.$\\
\noindent
$P_1$=$\{1,10,19,28,37,46\}$; $P_{2}$=$\{2,11,20,29,38,47\}$; $P_3$=$\{3,12,21,30,39,48\}$;\\
\noindent
$P_4$=$\{4,13,22,31,40,49\}$; $P_{5}$=$\{5,14,23,32,41,50\}$; $P_6$=$\{6,15,24,33,42,51\}$;\\
\noindent
$P_{7}$=$\{7,16,25,34,43,52\}$; $P_8$=$\{8,17,26,35,44,53\};$ $P_{9}$=$\{9,18,27,36,45\}$.\\
\noindent
$(iii)$~
$C(87,20)$: $a=4, r=7, k_1=13, k_2=6, r_2=1.$\\
\noindent
$P_1$=$\{1,68,48,28,8,75,55,35,15,82,62,42,22\}$; $P_2$=$\{2,69,49,29,9,76,56,36,16,83,63,43,23\}$;\\
\noindent
$P_3$=$\{3,70,50,30,10,77,57,37,17,84,64,44,24\}$; $P_4$=$\{4,71,51,31,11,78,58,38,18,85,65,45,25\}$;\\
\noindent
$P_5$=$\{5,72,52,32,12,79,59,39,19,86,66,46,26\}$; $P_6$=$\{6,73,53,33,13,80,60,40,20,87,67,47,27\}$; $P_7$=$\{7,74,54,34,14,81,61,41,21\}$.\\
\noindent
$(iv)$~$C(77,10)$: $a=7, r=7, k_1=3, r_1=1.$\\
\noindent
$P_1$=$\{1,68,58,48,38,28,18,8,75,65,55,45,35,25,15,5,72,62,52,42,32,22,12\}$;\\
\noindent
$P_2$=$\{2,69,59,49,39,29,19,9,76,66,56,46,36,26,16,6,73,63,53,43,33,23,13\}$;\\
\noindent
$P_3$=$\{3,70,60,50,40,30,20,10,77,67,57,47,37,27,17,7,74,64,54,44,34,24,14\}$;\\
\noindent
$P_4$=$\{4,71,61,51,41,31,21,11\}$.

\noindent
$(v)$~$C(56,9)$: $a=6, r=2, k_1=7, k_2=5, k_3=3, k_4=1, r_4=0.$\\
\noindent
$P_{1}$=$\{1,10,19,28,37,46,55,8,17,26,35,44,53,6,15,24,33,42,51,4,13,22,31,40,49\}$;\\
\noindent
$P_{2}$=$\{2,11,20,29,38,47,56,9,18,27,36,45,54,7,16,25,34,43,52,5,14,23,32,41,50\}$;\\
\noindent
$P_{3}$=$\{3,12,21,30,39,48\}.$

\section{The classification of circulant graphs $C(\mathbb{Z}_{n},\{k_1,k_2\})$}
In this section,  we mainly give the classification of circulant graphs $C(\mathbb{Z}_{n},\{k_1,k_2\})$ and discuss the dispersability of some cases.

\noindent \textbf{Definition 3.1.~}$^{[15]}~$\emph{Two graph $G,H$ are isomorphic $G\cong H$, if there are bijections $\theta:V(G)\rightarrow V(H)$ and $\phi:E(G)\rightarrow E(H)$ such that $\psi_{G}(e)=(u,v)$ if and only if $\psi_{H}(\phi(e))=(\theta(u),\theta(v))$, where $\psi_{G},\psi_{H}$ are the incidence functions for $G$ and $H$.}\\
\indent
For circulant graphs $G=C(\mathbb{Z}_{n},\{k_1,k_2\})$, before discussing the classification of $G$, we introduce a result about the Diophantine equation which is the main technique to judge which circle a vertex belongs to.

\noindent \textbf{Lemma 3.1.~}\emph{Given $a,b,c\in \mathbb{Z}^{+}$, if $gcd(a,b)=1$, then there exists a unique integer solution $(x_0,y_0)$ for the linear Diophantine equation $ax+by=c$ such that $0\leq x_0\leq b-1$. Furthermore, if $0\leq c \leq b-1$ and let $p(i)$ be the  ordinal number which represents the position  of the element $i$ in the ordered set $V_1^{\prime}=\{1,1+a,1+2a,\cdots,1+(b-1)a\}\pmod{b}$, then $p(1+c)$ is $1+x_0$.}\\
\noindent \textbf{Proof.~}Since $gcd(a,b)=1$, there exist integers $a^{\prime},b^{\prime}$ such that $a^{\prime}a+b^{\prime}b=1$, which is $(a^{\prime}c)a+(b^{\prime}c)b=c$. The integers $a^{\prime},b^{\prime}$ can be obtained by Euclidean algorithm and the general solution of the equation is $x=a^{\prime}c+tb,y=b^{\prime}c-ta(t \in \mathbb{Z})$. It is clear that there exists a unique solution $(x_0,y_0)$ such that $0\leq x_0\leq b-1$.~Since $gcd(a,b)=1$, $V_1^{\prime}$ is a rearrangement of $\mathbb{Z}_{b}$.~Assume $p(1+c)$ is $z_1$, then $1+c=1+a(z_1-1)-z_2b(z_2 \in \mathbb{Z}).$ That is $a(z_1-1)-z_2b=c.$ Hence $z_1=x_0+1.$
\hfill{$\square$}

For convenience, let $d_i=gcd(n,k_i)(i=1,2)$ and $d=gcd(d_1,d_2)$. We discuss the classification of $C(\mathbb{Z}_{n},\{k_1,k_2\})$ in the following four cases.

\vspace{0.2em}
\noindent\textbf{Case 1: $d_1=1$ or $d_2=1$}\\
\noindent \textbf{Lemma 3.2.~}\emph{If~$d_1=1$ or $d_2=1$, then~there~exists~an~integer~$k\in\mathbb{Z}_{n}$~such~that~$C(\mathbb{Z}_{n},\{k_1,k_2\})\cong C(n,k)$.}\\
\noindent \textbf{Proof.~}Without loss of generality,~we can assume $d_1=1$.~Let $k=x_0$~and~$(x_0,y_0)$ is the solution of the linear Diophantine~equation~$k_1x-ny=k_2$~such that~$0\leq x_0\leq n-1$.~Let $\theta:\mathbb{Z}_{n}\rightarrow \mathbb{Z}_{n}$ defined by $\theta(i)=p(i)(i\in\mathbb{Z}_{n})$ and $\phi:E(C(\mathbb{Z}_{n},\{k_1,k_2\}))$$\rightarrow E(C(n,k))$ defined by $\phi((i,j))=(\theta(i),\theta(j))(1 \leq i<j \leq n)$.~Because $d_1=1$, it is easy to see $\theta:\mathbb{Z}_{n}\rightarrow \mathbb{Z}_{n}$ is a bijection. Next we prove $\phi:E(C(\mathbb{Z}_{n},\{k_1,k_2\}))\rightarrow E(C(n,k))$ is well defined.\\
\indent
For $\forall (i,j)\in E(C(\mathbb{Z}_{n},\{k_1,k_2\}))$, we have $i-j\in \{\pm k_1,\pm k_2\}(\text{mod}~n)$. If $i-j\in\{\pm k_1\}(\text{mod}~n)$, by the definition of $\theta$, we have $\theta(i)-\theta(j)\in\{\pm1\}(\text{mod~}n)$. So $\phi((i,j))\in E(C(n,k)).$ If $i-j$ $\in\{\pm k_2\}(\text{mod}~n)$, we have $i=1+k_1(\theta(i)-1)-nt_1,j=1+k_1(\theta(j)-1)-nt_2(t_1,t_2 \in \mathbb{Z})$. Then,

\vspace{-1.4em}
$$i-j=k_1(\theta(i)-\theta(j))-n(t_1-t_2)\in\{\pm k_2\}(\text{mod~}n).$$
\vspace{-2em}

\noindent
By Lemma 3.1, $\theta(i)-\theta(j)\in \{\pm k\}(\text{mod~}n)$. Hence $\phi((i,j))\in E(C(n,k)).$ Because $\theta$ is a bijection, it is easy to see the mapping $\phi$ is both an injection and a surjection. The result is established.
\hfill{$\square$}

\vspace{0.2em}
\noindent\textbf{Case 2:  $d_1=d_2\neq 1$}\\
\noindent\textbf{Lemma 3.3.~}$^{[14]}$~\emph{The circulant graph $C(\mathbb{Z}_{n},\{k_1,k_2\})$ has $gcd(k_1,k_2,n)$ isomorphic connected components.}\\
\indent According to Lemma 3.2 and Lemma 3.3, we have the following result:\\
\noindent \textbf{Corollary 3.1.~}\emph{Let $n_1=n/d_1$.~The circulant graph $C(\mathbb{Z}_{n},\{k_1,k_2\})$ has $d$ isomorphic connected components.~And if $d_1=d_2$,~then each one of these connected components is isomorphic to $C(\mathbb{Z}_{n_1},\{k_1,k_2\})(\text{mod~}n_1)\cong C(n_1,k^{\prime})$, where $k^{\prime}=x^{\prime}$ and $(x^{\prime},y^{\prime})$ is the solution of the Diophantine equation $k_1x-n_1y=k_2$ such that $0\leq x^{\prime}\leq n_1-1$.}

\vspace{0.2em}
\noindent\textbf{Case 3: $1< d_1 < d_2<n/2$}

\noindent \textbf{Lemma 3.4.~}\emph{If $1 < d_1 < d_2<n/2$,~then the circulant graph $C(\mathbb{Z}_{n},\{k_1,k_2\})$ has $d$ isomorphic connected components and each one of these connected components is isomorphic to a Cartesian bundle~$C_{d_1/d}\Box^{\varphi}C_{n/d_1}$~over base~$C_{d_1/d}$ and fiber~$C_{n/d_1}$, with $\varphi$ as a cyclic~$x_1$-shift, where $(x_1,y_1)$ is the unique solution of the linear Diophantine~equation~$k_1x-(n/d)y=k_2d_1/d$~such that~$0\leq x_1\leq n/d_1-1$. Especially, $\varphi$ is trivial if and only if $d_1d_2\equiv 0\pmod n$.}

\noindent\textbf{Proof.~}(1)~If $d=1$, we provide an algorithm to show $C(\mathbb{Z}_{n},\{k_1,k_2\})$ is a Cartesian bundle $C_{d_1}\Box^{\varphi}C_{n/d_1}$.\\
\noindent \textbf{Step 1.~Find~$d_1$~cycles~induced~by~the~circulant~graph~$C(\mathbb{Z}_{n},\{k_1\})$.}
\vspace{-1em}
$$C_1:1\rightarrow 1+k_1\rightarrow 1+2k_1\rightarrow\cdots\rightarrow 1+(n/d_1-1)k_1~(\text{mod~}n),$$
\vspace{-2.3em}
\vspace{-0.4em}
$$C_2:2\rightarrow 2+k_1\rightarrow 2+2k_1\rightarrow\cdots\rightarrow 2+(n/d_1-1)k_1~(\text{mod~}n),$$
\vspace{-1.8em}
$$\vdots$$
\vspace{-2em}
$$C_{d_1}:d_1\rightarrow d_1+k_1\rightarrow d_1+2k_1\rightarrow \cdots\rightarrow d_1+(n/d_1-1)k_1~(\text{mod~}n).$$
\vspace{-2.1em}

\noindent \textbf{Step 2.~Rotate and rearrange the above $d_1$ cycles.}\\
\indent Let $A=\{1+ik_2~|~i\in\{1,2,3,\cdots,d_1-1,d_1\}\}\pmod{n}$, $B=\{1,2,3,\cdots,d_1-1,d_1\}$. In order to show the validity of the operation in this step, firstly, we construct a mapping $\phi$. Let $\phi:A\rightarrow B$ defined by $\phi(1+i k_2)=m$, where $m$ satisfies that $C_m$ includes the vertex $1+ik_2$, then the mapping $\phi:A\rightarrow B$ is a bijection and $\phi(1+d_1k_2)=1$. In fact, assume there are two integers $i,j(1\leq i<j \leq d_1)$ such that $\phi(1+i k_2)=\phi(1+j k_2)$, then the equation $1+ik_2+k_1x-ny=1+jk_2$ has solutions. That implies the following equation

\vspace{-1.6em}
\begin{equation*}\tag {*}
k_1x-ny=(j-i)k_2
\end{equation*}

\vspace{-0.4em}
\noindent has solution. Since~$1\leq j-i\leq d_1-1$,~$d_1\nmid j-i$.~It is easy to see $gcd(k_2,d_1)=1$,~so~$d_1\nmid (j-i)k_2$.~The linear Diophantine~equation $(*)$ has no solution. This is a contradiction.~So the mapping $\phi$ is an injection.~Furthermore, since $|A|=|B|$, $\phi$ is a surjection.~In addition,~it is easy to see that  there exist two integers $x,y$ such that $1+k_1x-ny=1+d_1k_2$, where $0\leq x\leq n/d_1-1$. Hence $\phi(1+d_1k_2)=1$.

Next, we give the methods for rotation and rearrangement.\\
\noindent
$(i)$~Place $C_1$ in the first line and place the cycle containing the vertex $1+k_2(i-1)(2\leq i \leq d_1)$ in the $ist$ line.\\
\noindent
$(ii)$~Keep the adjacency of the cycle $C_i(1\leq i \leq d_1)$ unchanged and rotate the cycle $C_i(2\leq i \leq d_1)$ until the vertex $1+k_2(i-1)(2\leq i \leq d_1)$ in the first column.

\noindent \textbf{Step 3.~Calculate the shift $\varphi$ of the Cartesian graph bundle of $C_{d_1}\Box^{\varphi} C_{n/d_1}$.}\\
\indent
After the above two steps,~by the definition of Cartesian graph bundle,~it is easy to see the circulant graph~$C(\mathbb{Z}_{n},\{k_1,k_2\})$~is a Cartesian graph bundle~$C_{d_1}\Box^{\varphi} C_{n/d_1}$,~where $\varphi$ is an~$x_1$-shift $(x_1 \in \mathbb{Z}$).~Next we calculate the integer~$x_1$. Because the vertex $1+d_1k_2$ belongs to $C_1$, it is sufficient to calculate the position difference of the vertices $1$ and $1+d_1k_2$ in $C_1$. Let $t_i=k_i/d_i(i=1,2)$. That is to find the solution of the linear Diophantine~equation

\vspace{-1.6em}
\begin{equation*}\tag {**}
k_1x-ny=d_1k_2,
\end{equation*}

\vspace{-0.6em}
\noindent which is reduced to $t_1x-(n/d_1)y=k_2$. It is clear that this equation exists a unique solution $(x_1,y_1)$ such that $0 \leq x_1 \leq n/d_1-1$. Thus $\varphi$ is a $x_1$-shift.

Finally, if $\varphi$ is trivial, obviously, $(0,y_1)$ is a solution of $(**)$. So we have $t_2d_1d_2\equiv0(\text{mod}~n)$. It is easy to see $gcd(t_2,n/(d_1d_2))=1$, so $d_1d_2\equiv 0(\text{mod}~n)$. If $d_1d_2\equiv 0(\text{mod}~n)$, by the equation $(**)$, we have $t_1d_1x\equiv 0(\text{mod}~n)$. Since $gcd(t_1,n/d_1)=1$, $x\equiv 0(\text{mod}~n/d_1)$. Thus $\varphi$ is trivial.

(2)~If $d\neq 1$, by Corollary 3.1, the circulant graph $C(\mathbb{Z}_{n},\{k_1,k_2\})$ has $d$ isomorphic connected components and the proof is similar to the case for $d=1$.
\hfill{$\square$}

\noindent \textbf{Example.~}$C(\mathbb{Z}_{60},\{28,35\})\cong C_4\Box^{\varphi} C_{15}$, where $\varphi$ is a $5$-shift(see Fig.3).

\noindent \textbf{Step 1.~}
\vspace{-1em}
\noindent$$ C_1:01\rightarrow 29\rightarrow 57 \rightarrow 25\rightarrow 53 \rightarrow 21 \rightarrow 49 \rightarrow 17 \rightarrow 45 \rightarrow 13 \rightarrow41 \rightarrow 09 \rightarrow 37 \rightarrow 05 \rightarrow33;$$
\vspace{-2.4em}
\noindent$$ C_2:02\rightarrow 30\rightarrow 58 \rightarrow 26\rightarrow 54 \rightarrow 22 \rightarrow 50 \rightarrow 18 \rightarrow 46 \rightarrow 14 \rightarrow42 \rightarrow 10 \rightarrow 38 \rightarrow 06 \rightarrow34;$$
\vspace{-2.1em}
\noindent$$ C_3:03\rightarrow 31\rightarrow 59\rightarrow 27\rightarrow 55 \rightarrow 23 \rightarrow 51 \rightarrow 19 \rightarrow 47 \rightarrow 15 \rightarrow43 \rightarrow 11 \rightarrow 39 \rightarrow 07 \rightarrow35;$$
\vspace{-2.1em}
\noindent$$ C_4:04\rightarrow32\rightarrow 60\rightarrow 28\rightarrow 56 \rightarrow 24 \rightarrow 52 \rightarrow 20 \rightarrow 48 \rightarrow 16 \rightarrow44 \rightarrow 12 \rightarrow 40 \rightarrow 08 \rightarrow36.$$
\vspace{-2em}

\noindent \textbf{Step 2.~}
\vspace{-0.8em}
\noindent$$ C_1:01\rightarrow 29\rightarrow 57 \rightarrow 25\rightarrow 53 \rightarrow 21 \rightarrow 49 \rightarrow 17 \rightarrow 45 \rightarrow 13 \rightarrow41 \rightarrow 09 \rightarrow 37 \rightarrow 05 \rightarrow33;$$
\vspace{-2.4em}
\noindent$$ C_4:36 \rightarrow 04\rightarrow32\rightarrow 60\rightarrow 28\rightarrow 56 \rightarrow 24 \rightarrow 52 \rightarrow 20 \rightarrow 48 \rightarrow 16 \rightarrow44 \rightarrow 12 \rightarrow 40 \rightarrow 08;$$
\vspace{-2.1em}
\noindent$$ C_3: 11 \rightarrow 39 \rightarrow 07 \rightarrow35 \rightarrow 03\rightarrow 31\rightarrow 59\rightarrow 27\rightarrow 55 \rightarrow 23 \rightarrow 51 \rightarrow 19 \rightarrow 47 \rightarrow 15 \rightarrow43;$$
\vspace{-2.1em}
\noindent$$ C_2:46 \rightarrow 14 \rightarrow42 \rightarrow 10 \rightarrow 38 \rightarrow 06 \rightarrow34 \rightarrow  02\rightarrow 30\rightarrow 58 \rightarrow 26\rightarrow 54 \rightarrow 22 \rightarrow 50 \rightarrow 18.$$

\vspace{-1em}
\begin{figure}[htbp]
\centering
\includegraphics[height=3.3cm, width=0.55\textwidth]{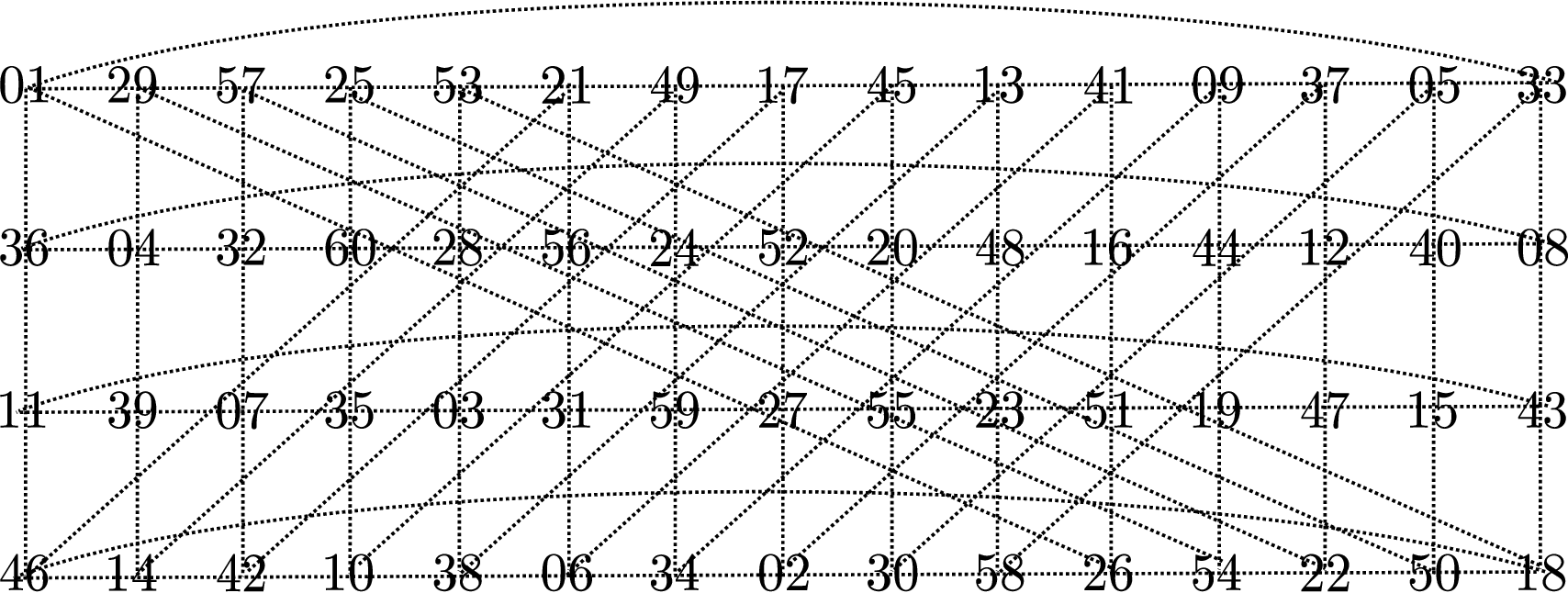}
\centerline{Fig.3  ~The circulant graphs $C(\mathbb{Z}_{60},\{28,35\})\cong C_4\Box^{\varphi} C_{15}$, where $\varphi$ is a $5$-shift.}%
\end{figure}
\vspace{-1.2em}

\noindent \textbf{Step 3.~}$x=5,y=0$~is the unique solution of the equation $7x-15y=35$ such that $0\leq x\leq 14$.

\vspace{0.2em}
\noindent\textbf{Case 4: $n$ is even and $1< d_1 < d_2=n/2$}\\
\noindent \textbf{Lemma 3.5.~}\emph{Let $n_1=n/d_1$. If $n$ is even and $1 < d_1 < d_2=n/2$, then}\\
\noindent
\emph{(i)~If $n_1$ is odd, then $C(\mathbb{Z}_{n},\{k_1,k_2\})$ has $d=gcd(d_1,d_2)=d_1/2$ isomorphic connected components, where each connected component is isomorphic to the Cartesian product $K_2\Box C_{n_1}$ of the complete graph $K_2$ and the cycle $C_{n_1}$;}\\
\noindent
\emph{(ii)~If $n_1$ is even, then $C(\mathbb{Z}_{n},\{k_1,k_2\})$ has $d=gcd(d_1,d_2)=d_1$ isomorphic connected components, where each connected component is isomorphic to the circulant graph $C(n_1,n_1/2)$.}\\
\noindent \textbf{Proof.~}$(i)$~Because $n$ is even and $n_1$ is odd,~$d_1$ is even and $d_1\nmid d_2$.~Since $d_1|n$, $d_1/2|n/2$.~Therefore $d=gcd(d_1,d_2)=d_1/2$.~By Corollary 3.1,~the circulant graph $C(\mathbb{Z}_{n},\{k_1,k_2\})$ has $d=d_1/2$ isomorphic connected components.~We can obtain the cycles $C_1,C_2,\cdots,C_d$ induced by $C(\mathbb{Z}_{n},\{k_1\})$ by the proof of Lemma 3.4.
Since $gcd(k_1,n_1)=1$, by Lemma 3.1, there exists a unique solution $(x_0,y_0)$ for the equation $i+n/2=i+d_1/2+k_1x-n_1y$, where $0\leq x_0\leq n_1-1$.  Hence for each vertex $i(1\leq i \leq d_1/2)$, its adjacent vertex $i+n/2$ belongs to the cycle $C_{i+d_1/2}$. It is easy to see $C_i$ and $C_{i+d_1/2}(1\leq i\leq d_1/2)$ are connected by a perfect matching and each connected component is isomorphic to $K_2\Box C_{n_1}$.\\
\noindent $(ii)$~Since $n/d_1$ is even,~$2d_1\mid n$,~which implies that $d_1\mid n/2$.~So we have $d=d_1$.~So  $C(\mathbb{Z}_{n},\{k_1,k_2\})$ has $d_1$ isomorphic connected components. In order to show each connected component is isomorphic to $C(n_1,n_1/2)$, it suffices to show $1+n/2$ belongs to $C_1$, and the distance between the vertices $1$ and $1+n/2$ is $n_1/2$. Since $gcd(k_1,n_1)=1$, there exist two integers $x_0,y_0$ such that $1+k_1x-n_1y=1+n/2$ holds, where $0\leq x_0\leq n_1-1$. In fact, it is easy to see $x_0=n_1/2$, $y_0=(k_1-d_1)/2$ is the unique integer solution. The result is established.
\hfill{$\square$}

\vspace{0.4em}
In summary, we obtain the following classification of circulant graphs $C(\mathbb{Z}_{n},\{k_1,k_2\})$.

\noindent \textbf{Theorem 3.1.~}\emph{For the circulant graph $C(\mathbb{Z}_{n},\{k_1,k_2\})$, we have}
\vspace{-0.4em}
$$C(\mathbb{Z}_{n},\{k_1,k_2\})\cong\left\{\begin{array}{l}
\underset{i=1}{\overset{gcd(d_1,d_2)}{\cup}}C(n/d_{1},k), \ \ \ \ \ \ \ \ \ \ \ \ \ \ \ \ \ \ \ \ d_1=1~\text{or}~d_1=d_2\neq1,\\
\underset{i=1}{\overset{gcd(d_1,d_2)}{\cup}}C_{d_1/gcd(d_1,d_2)}\Box^{\varphi} C_{n/d_1}, \ \ \ \ \ \ \ \ 1<d_1<d_2<n/2,\\
\underset{i=1}{\overset{gcd(d_1,d_2)}{\cup}}K_2\Box C_{n/d_1},\ \ \ \ \ \ \ \ \ \ \ \ \ \ \ \ \ \ \ \ \ n \text{~is even},1<d_1<d_2=n/2,2\nmid (n/d_1),\\
\underset{i=1}{\overset{gcd(d_1,d_2)}{\cup}} C(n/d_1,n/(2d_1)),\ \ \ \ \ \ \ \ \ \ \ \ n \text{~is even},1<d_1<d_2=n/2,2\mid (n/d_1),\\
\end{array} \right.$$

\vspace{-0.3em}
\noindent
where~$k=x_0$~and~$(x_0,y_0)$ is a solution of the linear Diophantine~equation~$k_1x-(n/d_1)y=k_2$~such that~$0\leq x_0\leq n/d_1-1$; 
$\varphi$ is a cyclic $x_1$-shift and~$(x_1,y_1)$ is a solution of the linear Diophantine~equation~$k_1x-(n/gcd(d_1,d_2))y=k_2d_1/gcd(d_1,d_2)$~such that~$0\leq x_1\leq n/d_1-1$. Especially, $\varphi$ is trivial if and only if $d_1d_2\equiv 0(\text{mod}~n)$.\\
\noindent \textbf{Remark 2.~}The matching book thickness of the Cartesian product of two cycles has been solved in~$[4],[8],[9],[10]$,~the matching book embedding of~$K_2\Box C_t$ has been solved in $[16]$.~The dispersability of the circulant graph $C(n,k)$ will be discuss in the following sections.

\section{The case: $n$ is even}
By Lemma 2.6, the circulant graph $C(n,\{k\})$ is the disjoint union of $C_1,C_2,\cdots,C_d$. Let $V_i$ be the clockwise cyclic ordered vertex set of the cycle $C_i$, specifically, $V_i=\{i, i+k, i+2k, \cdots, i+(n/d-1)k\}(1 \leq i \leq d)$. According to the parity of $k$, we consider two cases to discuss the dispersability of the circulant graph $C(n,k)$.

\vspace{0.4em}
\noindent \textbf{Theorem 4.1.~}\emph{Let $G=C(n,k)$, where $n$ is even, $k$ is odd, then $mbt(G)=\Delta(G).$}\\
\noindent \textbf{Proof.}
By Lemma 2.5, it suffices to show that $mbt(G)\leq\Delta(G)$.~Put the vertices of $G$ clockwise along a circle in the order $1,n,3,n-2,\cdots,n-3,4,n-1,2$.

\vspace{-0.5em}
\begin{figure}[htbp]
\centering
\includegraphics[height=3.6cm, width=0.62\textwidth]{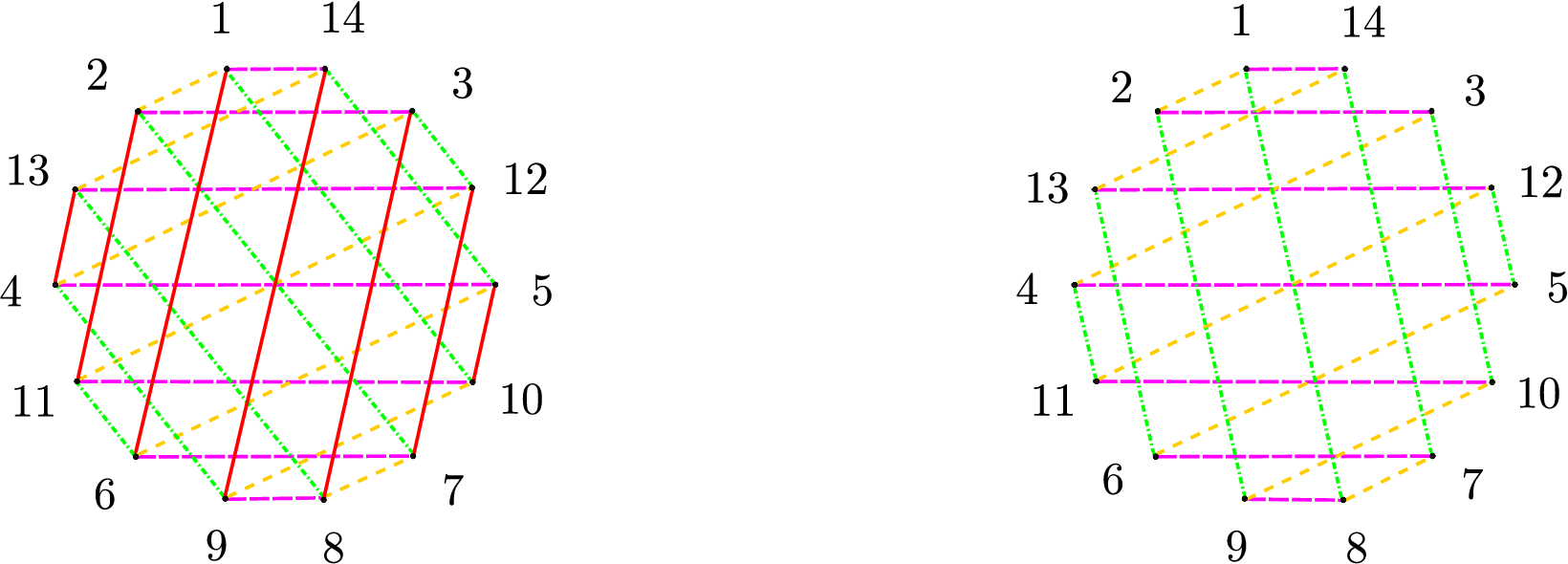}
\centerline{Fig.4  ~The matching book embedding of $C(14,5)$(left) and $C(14,7)$(right).}%
\end{figure}
\vspace{-1.2em}

If $k \neq n/2$, we have $\Delta(G)=4$.~The edges of $G$ can be colored well with $\Delta(G)$ colors in the following way(see Figure $4$ for $C(14,5)$(left) and $C(14,7)$(right)):

Color the edges of $\{(i,i+1)~|~1 \leq i \leq n, i \text{~is~odd}\}$ yellow,  the edges of $\{(i,i+1)~|~1 \leq i \leq n, i \text{~is~even}\}$ purple, the edges of $\{(i,i+k)~|~1 \leq i \leq n, i \text{~is~even}\}$ green, the edges of $\{(i,i+k)~|~1 \leq i \leq n, i \text{~is~odd}\}$ red.

If $k=n/2$, we have $\Delta(G)=3$.~For each vertex of $G$, the edges colored with green and red are coincide. Hence $mbt(G)\leq\Delta(G)$.
\hfill{$\square$}

\noindent \textbf{Theorem 4.2.~}\emph{Let $G=C(n,k)$, where $n$ and $k$ are both even, then $mbt(G)=\Delta(G)+1.$}

\noindent \textbf{Proof.}
Since $n,k$ are both even, $d=gcd(n,k)$ is even. According to Lemma 2.5, it is sufficient to show that $mbt(G)\leq\Delta(G)+1$ .

\noindent\textbf{Case 1: $k=n/2$}

Put the vertices of $G$ counterclockwise along a circle in the order $1,2,\cdots,n/2,n,n-1,n-2,\cdots,n/2+1$.~All edges can be colored well in the following way(see Fig.5 for $C(12,6)$):

\indent
Color the edge  $(1,n)$ purple, the edges of $\{(i,i+1)~|~i \in \{1,3,\cdots,n-1\}\}$ red, the edges of $\{(i,i+1)~|~i \in\{2,4,\cdots,n-2\}\}$ green, the edges of $\{(i,i+k)~|~1\leq i \leq n/2\}$ yellow.

\vspace{-0.4em}
\begin{figure}[htbp]
\centering
\includegraphics[height=2.6cm, width=0.26\textwidth]{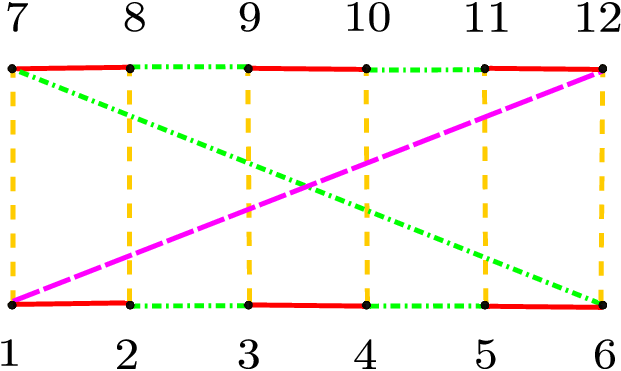}
\centerline{Fig.5  ~The matching book embedding of $C(12,6)$.}%
\end{figure}
\vspace{-1em}

\noindent\textbf{Case 2: $k\neq n/2$}

Put the vertices of $G$ counterclockwise along a circle in the order $V_{1}^{-}V_{2}V_{3}^{-}V_{4}\cdots V_{d-1}^{-}V_{d}$. The edges of $G$ can be colored well with $\Delta(G)+1$ colors in the following two steps(see Fig.6 for $C(24,4)$(left) and $C(56,8)$(right)):
\vspace{-0.8em}
\begin{figure}[htbp]
\centering
\includegraphics[height=5cm, width=0.76\textwidth]{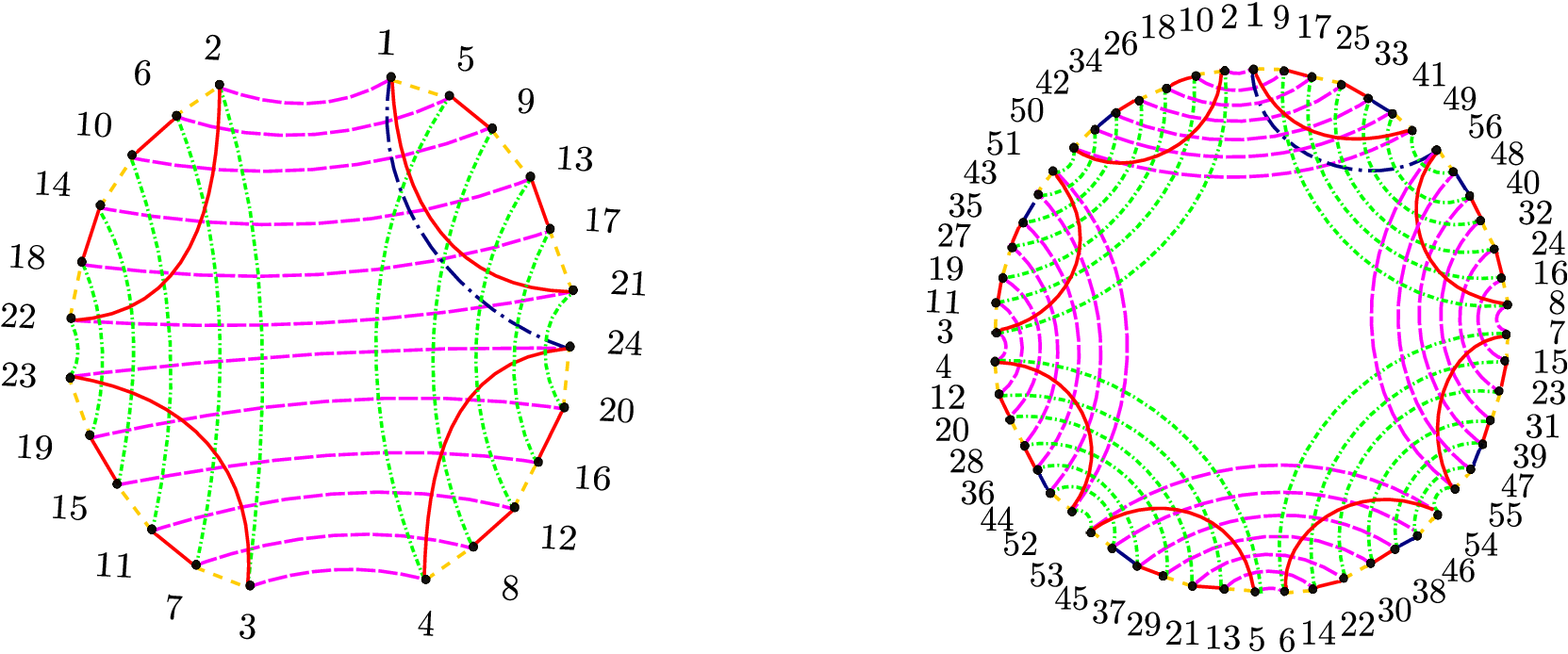}
\centerline{Fig.6  ~The matching book embedding of $C(24,4)$(left) and $C(56,8)$(right).}%
\end{figure}

\vspace{-1em}
\noindent
\textbf{Step 1:} The coloring of the $n$-cycle $C_0$.\\
\noindent
Purple: $\{(j,j+1)~|~j \in V_{i}, i \in\{1,3,\cdots,d-1\} \}$;\\
\noindent
Green:~$\{(j,j+1)~|~j\in V_{i}, i\in\{2,4,\cdots,d-2\} \}$, $\{(j,j+1)~|~j \text{~is before the element $n$ in}~V_{d}\}$;\\
\noindent
Blue: $\{(j,j+1)~|~j \text{~is not before the element $n$ in}~V_{d}\}$.\\
\noindent
\textbf{Step 2:} The coloring of the cycles $C_1,C_2,\cdots,C_d$.

Color the edges $(i,i+(n/d-1)k)(1 \leq i \leq d)$ red.~Regardless of the parity of $n/d$, it is easy to use red, blue and yellow to color the remaining edges induced by $V_i(1 \leq i \leq d)$.

Thus, $mbt(G)\leq\Delta(G)+1.$ The result is established.
\hfill{$\square$}


\section{The case: $n$ is odd, $d=gcd(n,k) \neq 1$}

\noindent \textbf{Theorem 5.1.~}\emph{Let $G=C(n,k)(k \leq \lfloor n/2 \rfloor)$, where $n$ is odd and $d \neq 1$, then $mbt(G)=5.$}

\noindent \textbf{Proof.}
According to Lemma 2.5, 
it is sufficient to show that $mbt(G)\leq\Delta(G)+1$.\\
\indent
Let the ordered vertex set $V_{d-1}^{\prime~}$ be a rearrangement of $V_{d-1}\cup V_d$, specifically, $V_{d-1}^{\prime}$= $\{d-1,d,d-1+k,d+k,d+2k,d-1+2k,d-1+3k,d+3k,d+4k,d-1+4k,\cdots,d-1+(n/d-2)k,d+(n/d-2)k,d+(n/d-1)k,d-1+(n/d-1)k\}$.

\vspace{0.4em}
\noindent\textbf{Case 1: $n=3k$}

Put the vertices of $G$ clockwise along a circle in the order $V_1V_2V_3^{-}V_4\cdots V_{d-2}^{-}V_{d-1}V_{d}^{-}$.~The edges of $G$ can be colored well with $\Delta(G)+1$ colors in the following way(see Fig.7 for $C(9,3)$(left) and $C(27,9)$(right)):

Yellow: $\{(j,j+1)~|~j \in V_i, i \in\{2,4,\cdots,d-1\}\}$, $\{(1,1+2k)\}$;\\
\indent
Red: $\{(j,j+1)~|~j \in V_i, i \in \{3,5,7,\cdots,d-2\}\}$, $\{(1,1+k),(2,2+k),(1+2k,2+2k),(d,d+2k)\}$;\\
\indent
Purple: $\{(1,2)\}$, $\{(i+k,i+2k)~|~i \in\{1,2,d\}\}$, $\{(i,i+k)~|~i \in \{3,4,\cdots,d-1\}\}$;\\
\indent
Blue: $\{(1,n),(1+k,2+k),(d,d+k)\}$, $\{(i,i+2k)~|~i \in \{3,4,\cdots,d-1\}\}$;\\
\indent
Green: $\{(1+k,k),(1+2k,2k),(2,2+2k)\}$, $\{(i+k,i+2k)~|~i \in \{3,4,\cdots,d-1\}\}$.

\vspace{-0.6em}
\begin{figure}[htbp]
\centering
\includegraphics[height=4.1cm, width=0.68\textwidth]{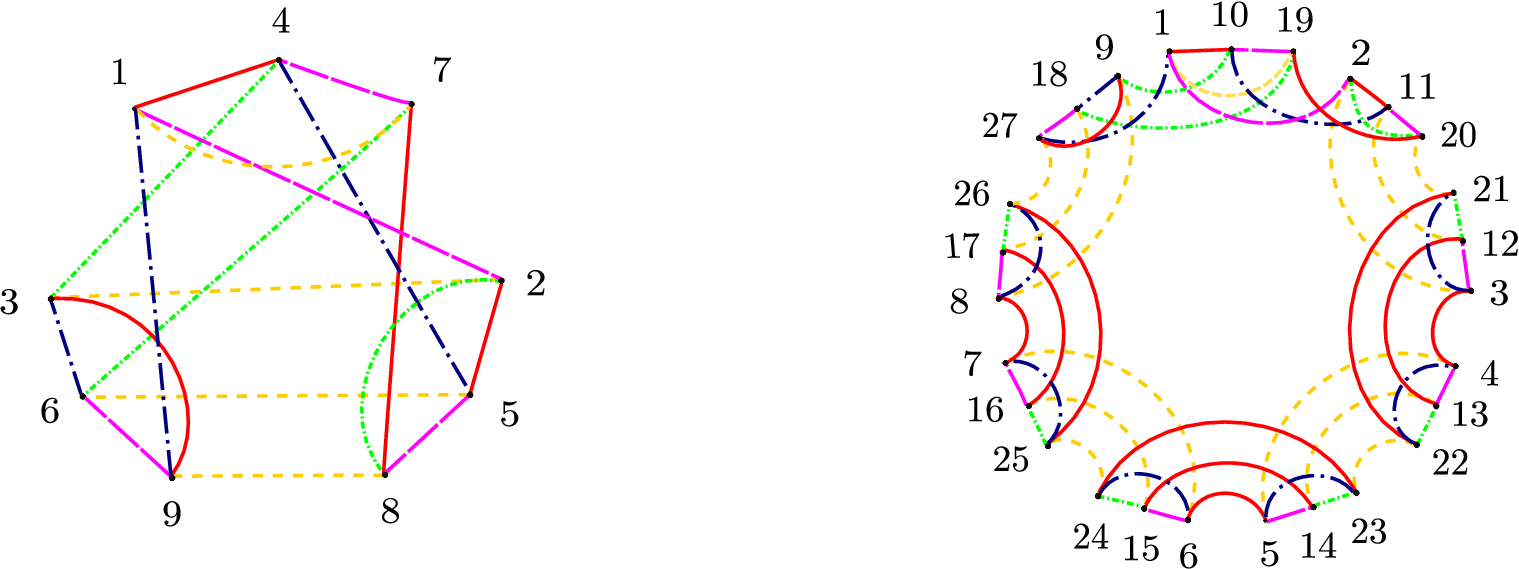}
\centerline{Fig.7 ~The matching book embedding of $C(9,3)$(left) and $C(27,9)$(right).}
\end{figure}
\vspace{-1em}

\vspace{0.4em}
\noindent\textbf{Case 2: $n\neq 3k,d=3$}

\indent
Put the vertices of $G$ clockwise along a circle in the order $V_1V_2^{{\prime}-}$. The edges of $G$ can be colored well with $\Delta(G)+1$ colors in the following four steps(see Fig.8 for $C(27,3)$(left) and $C(27,6)$(right)):\\
\noindent
\textbf{Step 1:~}The coloring of the edges between $C_1$ and $C_3$, the edges between $C_1$ and $C_2$.

Purple: $\{(i,i-1)~|~i \text{~is~the~element~of~$V_1$~before}~4\}$;

Green: $\{(i,i-1)~|~i \text{~is~the~element~of~$V_1$~after}~4+k\}$, $\{(1,2),(3,4)\}$;

Yellow: $\{(i,i+1)~|~i\in V_1\backslash\{1,1+k\}\}$;

Blue: $\{(1+k,2+k),(3+k,4+k)\}$.\\
\noindent
\textbf{Step 2:~}The coloring of the edges between $C_2$ and $C_3$.\\
\indent
Color the edge $(2,3)$ blue, the edge $(2+k,3+k)$ green, the edge $(2+(n/d-2)k,3+(n/d-2)k)$ red. As for the remaining edges between $C_2$ and $C_3$, if $d=k$, color the edge $(2+(n/d-3)k,3+(n/d-3)k)$ blue, the edge $(2+(n/d-1)k,3+(n/d-1)k)$ green, the remaining edges purple; otherwise, color the edge $(2+(n/d-1)k,3+(n/d-1)k)$ blue, the edges $(j-1,j)$ purple, where $j \in V_d$, $j$ is before $n$, and color other edges green.

\vspace{-0.8em}
\begin{figure}[htbp]
\centering
\includegraphics[height=4.2cm, width=0.68\textwidth]{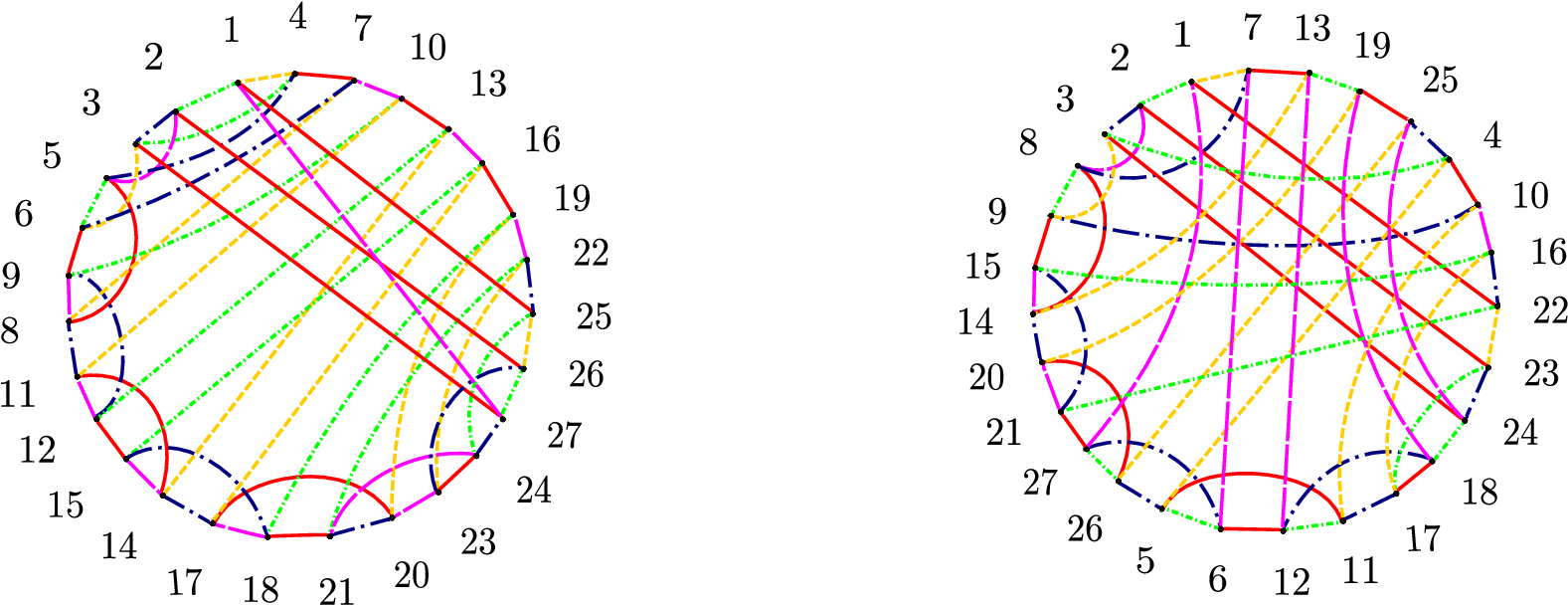}\\
\centerline{Fig.8 ~The matching book embedding of $C(27,3)$(left) and $C(27,6)$(right).}%
\end{figure}
\vspace{-1em}

\noindent
\textbf{Step 3: }The coloring of the cycles $C_2,C_3$.\\
\indent
Color the edge $(2,2+k)$ purple, the edge $(3,3+k)$ yellow, the edges of $\{(i-1,i-1+k),(i,i+k)~|~i=3+jk, j\in \{2,4,\cdots,n/d-5\}\}$ blue, the edges of $\{(i-1,i-1+k),(i,i+k)~|~i=3+jk, j\in \{1,3,\cdots,n/d-4\}\cup \{n/d-1\} \}$ red. If $d=k$, color the edges of $\{(i-1,i-1+k),(i,i+k)~|~i=3+jk, j=n/d-3\}$ purple, the edges of $\{(i-1,i-1+k),(i,i+k)~|~i=3+jk, j=n/d-2\}$ blue. Otherwise, color the edges of $\{(i-1,i-1+k),(i,i+k)~|~i=3+jk, j=n/d-3\}$ blue, the edges of $\{(i-1,i-1+k),(i,i+k)~|~i=3+jk, j=n/d-2\}$ green.

\vspace{-0.8em}
\begin{figure}[htbp]
\centering
\includegraphics[height=5cm, width=0.56\textwidth]{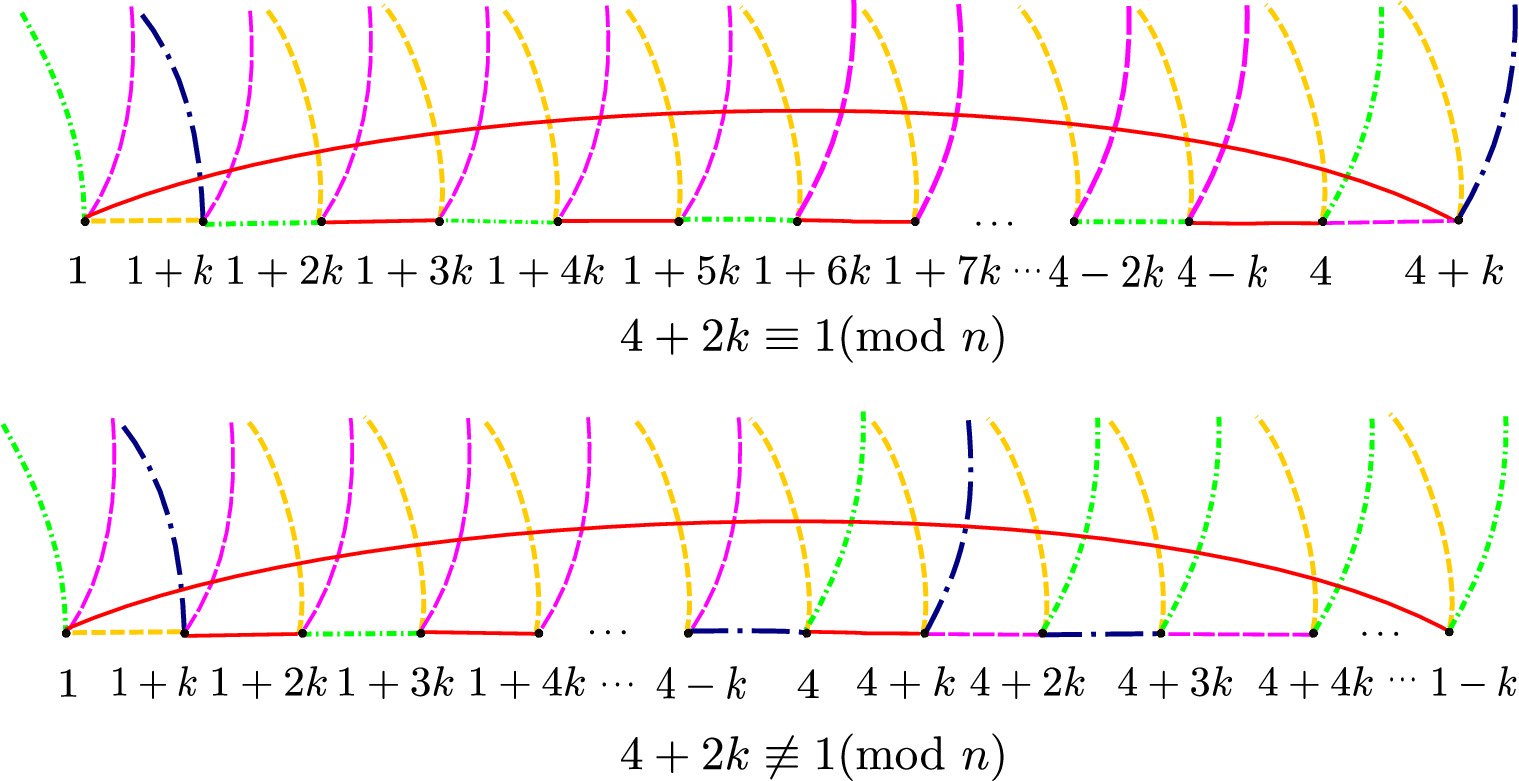}
\centerline{Fig.9  ~The edge-coloring of $C_1$.}
\end{figure}

\vspace{-1.1em}
\noindent
\textbf{Step 4: }The coloring of the cycle $C_1$.\\
\indent Color the edge $(1,1+k)$ yellow, the edge $(1,1-k)$ red. If $d=k$, color the edge $(n-2-k,n-2)$ blue,~the remaining edges of $C_1$ with red and purple alternately. Otherwise, the coloring of the edges incident to $V_1$ has the following characteristics (ignore the vertex $1$).
All vertices of $V_1$ before the vertex $4$ incident to a purple edge,~the vertex $4+k$ incident to a blue edge, the other vertices incident to a green edge. In addition, the vertices $1+ki(i\in \{2,3,\cdots,n/d-1\})$ incident to a yellow edge, the vertex $1+k$ has a blue edge(see Fig.9).\\
\indent If $4+2k\equiv 1(\text{mod}~n)$, color the edge $(4,4+k)$ purple and the remaining even edges induced by $V_1$ with red and green alternately, where the edge $(4-k,4)$ is red.~Otherwise, color the edge $(4,4+k)$ red, the edge $(4,4-k)$ blue, the edges induced by $V_1$ between the vertices $1+k,4-k$ with red and green alternately, the remaining edges with purple and blue alternately, where the edge $(4+k,4+2k)$ is purple.

\vspace{0.4em}
\noindent\textbf{Case 3: $n\neq 3k,d \geq 5$}\\
\indent
Put the vertices of $G$ counterclockwise along a circle in the order $V_1V_3\cdots V_{d-2}V_{d-1}^{{\prime}-}V_{d-3}^{-}\cdots V_4^{-}\\V_2^{-}$.~All edges can be colored well with $\Delta(G)+1$ colors in the following four steps(see Fig.10 for $C(65,5)$(left) and $C(65,15)$(right)):

\vspace{-0.9em}
\begin{figure}[htbp]
\centering
\includegraphics[height=5.6cm, width=0.96\textwidth]{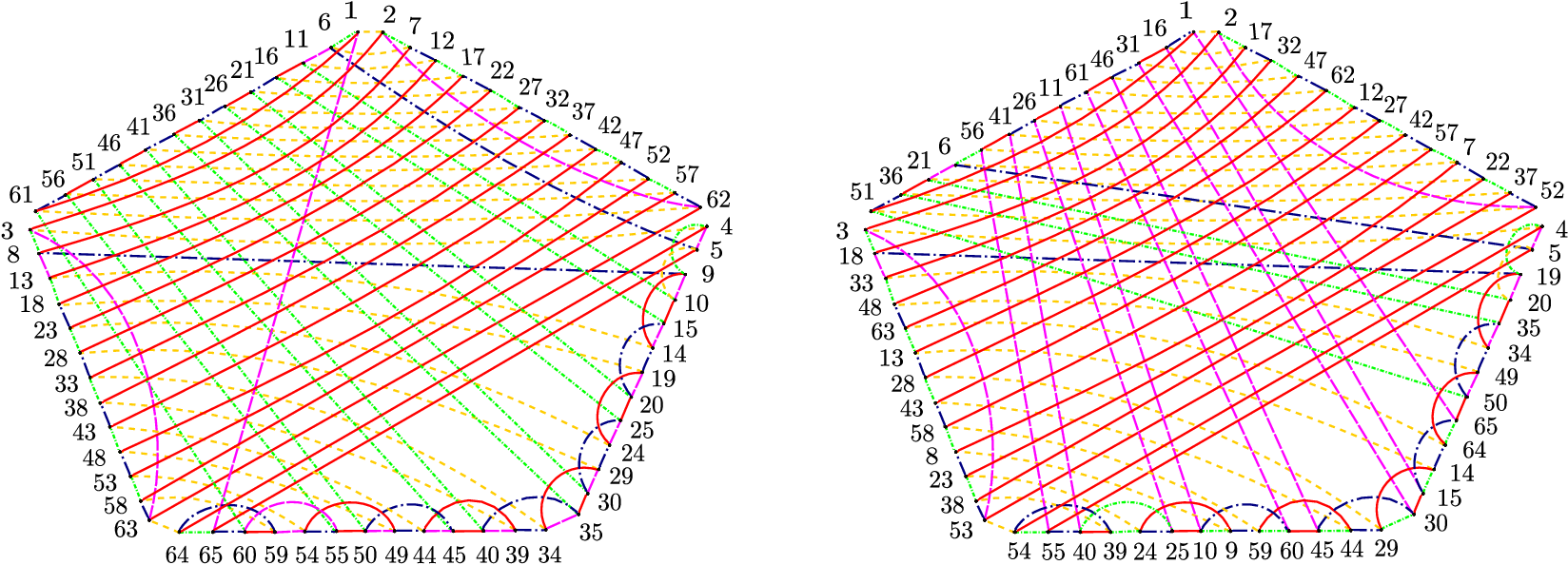}
\centerline{Fig.10 ~The matching book embedding of $C(65,5)$(left) and $C(65,15)$(right).}%
\end{figure}
\vspace{-1.2em}

\noindent
\textbf{Step 1: }The coloring of the edges between $C_{d-1}$ and $C_d$.\\
\indent
Color the edge $(d-1+(|V_d|-2)k,d+(|V_d|-2)k)$ red, the edge $(d-1+(|V_d|-3)k,d+(|V_d|-3)k)$ blue. As for the remaining edges, if $j$ is before the element $n$ in $V_d$, color the edges $(j-1,j)$ purple and color other edges green.\\
\noindent
\textbf{Step 2: }The coloring of the other edges of the $n$-cycle $C_0$.

Red: $\{(j,j+1)~|~j \in V_i, i \in \{2,4,\cdots,d-3\}\}$;

Yellow: $\{(j,j+1)~|~j \in V_i, i\in \{1,3,\cdots,d-2\}, j \neq d-2+k\}$;

Blue: $\{(d-2+k,d-1+k),(d,d+1)\}$;

Green: $\{(j,j+1)~|~\text{$j$ is before the element $n$ in $V_d$},j \neq d\}$;

Purple: $\{(j,j+1)~|~\text{$j$ is not before the element $n$ in $V_d$}\}$.

\noindent
\textbf{Step 3: }The coloring of the cycles $C_{d-1},C_d$.\\
\indent Color the edge $(d-1,d-1+k)$ green, the edge $(d,d+k)$ yellow, the edges of $\{(i-1,i-1+k),(i,i+k)~|~i=d+jk, j \in \{1,3,\cdots,|V_d|-4\}\cup\{|V_d|-1\}\}$ red, the edges of $\{(i-1,i-1+k),(i,i+k)~|~i=d+jk,j \in \{2,4,\cdots,|V_d|-5\}\cup\{|V_d|-2\}\}$ blue.
If $d+(|V_d|-3)k$ is before the element $n$ in $V_d$, color the edges of $\{(i-1,i-1+k),(i,i+k)~|~i=d+(|V_d|-3)k\}$ purple. Otherwise, color them green.\\
\noindent\textbf{Step 4: }The coloring of the cycles $C_1,C_2,\cdots,C_{d-2}$.\\
\indent As for the cycle $C_1$, color the edge $(d+1-k,d+1)$ green, the edge $(d+1,d+1+k)$ purple, and color the remaining edges with red and blue alternately, where the edge $(1,1-k)$ is red. Color the edges of $\{(i,i-k)~|~i \in \{2,3,\cdots,d-2\}\}$ purple. As for the edges of the cycle $C_j,j=2,3,\cdots,d-3$, color the remaining edges with blue and green alternately. Color the edge $(d-2,d-2+k)$ green, the edge $(d-2+k,d-2+2k)$ purple, and the remaining edges with green and blue alternately.

All edges of $G$ can be matching book embedded in five pages. Therefore, $mbt(G)=5$.
\hfill{$\square$}

\section{The case: $n$ is odd, $d=gcd(n,k)=1$}
Let $n=ak+r$, the vertices of $C_1$ can be divided into $P_1,P_2,\cdots,P_t$ as in Section 2. If $r=1$, it is easy to check that $|P_1|=a+1$, $|P_i|=a(2 \leq i \leq k)$.

\vspace{0.4em}
\noindent \textbf{Theorem 6.1.~}\emph{Let $G=C(n,k)$, where $n$ is odd. If~$r=1$, then $mbt(G)=\Delta(G)+1=5.$}

\noindent \textbf{Proof.}
By Lemma 2.5, it is sufficient to show that $mbt(G)\leq\Delta(G)+1$.

\vspace{0.4em}
\noindent\textbf{Case 1: $a=2$}

Let $Q_1,Q_2$ be ordered vertex sets, where $Q_1=\{1,2,\cdots,k+1\}$, $Q_2=\{k+2,k+3,\cdots,n\}.$ Put the vertices of $G$ counterclockwise along a circle in the order $Q_1Q_2^{-}$. The edges of $G$ can be colored well with $\Delta(G)+1$ colors in the following two steps(see Fig.11 for $C(27,13)$(left) and $C(25,12)$(right)):

\noindent
\textbf{Step 1:} The coloring of the $n$-cycle $C_1$.

Color the edge $(1,1+k)$ purple, the edges of $\{(i,i+k)~|~i \in Q_1\backslash\{1\} \}$ green, the edges of $\{(i,i+k)~|~i \in Q_2 \}$ red.

\noindent
\textbf{Step 2:} The coloring of the $n$-cycle $C_0$.\\
\indent
Color the edges $(1,2),(1+k,2+k)$ yellow, the edge $(1,n)$ blue, the path induced by $Q_1\backslash\{1\}$ and the path induced by $Q_2$ with blue and purple alternately, where the edge $(k,k+1)$ is blue, the edge $(n,n-1)$ is purple.

\vspace{-0.6em}
\begin{figure}[htbp]
\centering
\includegraphics[height=4.6cm, width=0.7\textwidth]{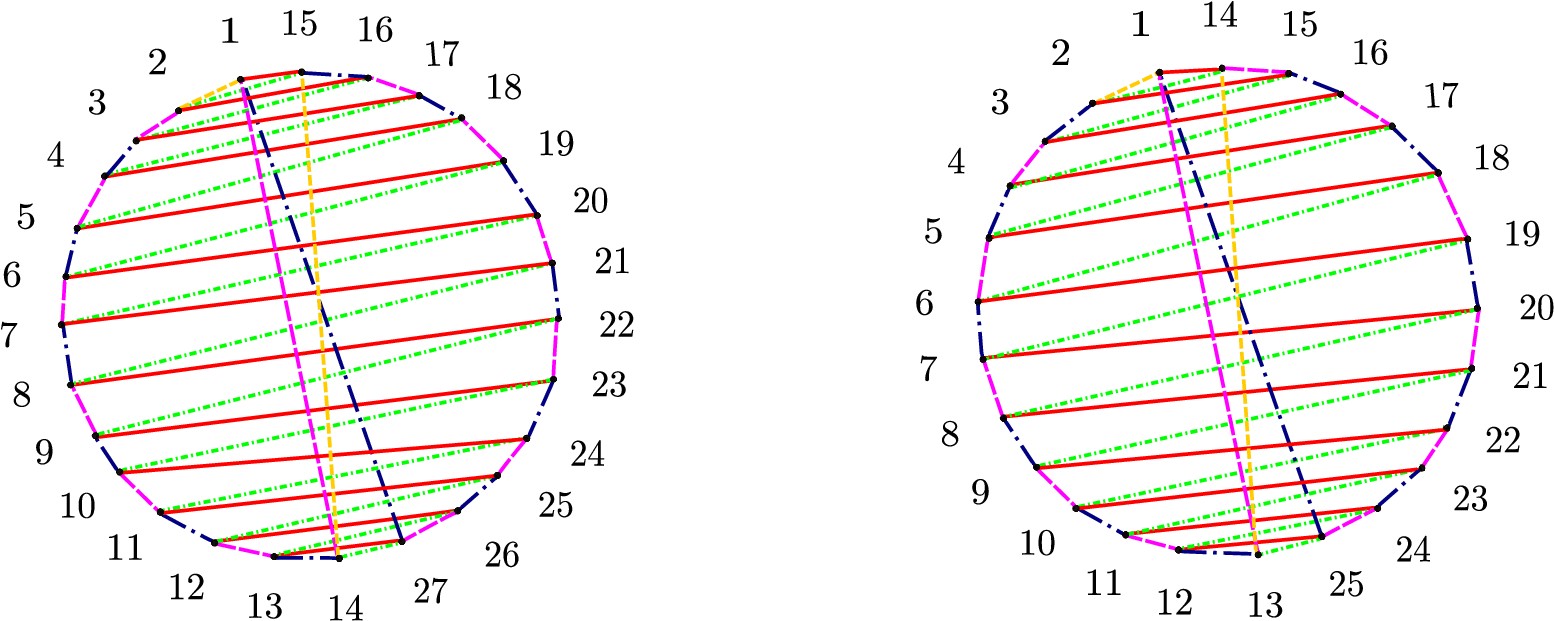}
\centerline{Fig.11  ~The matching book embedding of $C(27,13)$(left) and $C(25,12)$(right).}%
\end{figure}
\vspace{-1.2em}

\vspace{0.4em}
\noindent\textbf{Case 2: $a \geq 3$}\\
\vspace{0.3em}
\noindent\textbf{Subcase 2.1. $k$ is even.}\\
\indent Put the vertices of $G$ counterclockwise along a circle in the order $P_1^{-}P_2P_3^{-}\cdots P_{k}$.~The edges of $G$ can be colored well with $\Delta(G)+1$ colors in the following three steps(see Fig.12 for $C(25,8)$(left) and $C(25,6)$(right)):\\
\noindent
\textbf{Step 1:} The coloring of the $n$-cycle $C_0$.\\
\indent
Color the edge $(1,n)$ yellow, the edges of $\{(i,i-1)~|~i \in P_j, j\in\{2,4,\cdots,k\}\}$ green, the edges of $\{(i,i+1)~|~i \in P_j, j\in\{2,4,\cdots,k\}\}$ red.\\
\noindent
\textbf{Step 2:} The coloring of the remaining edges incident to the endpoints of the path induced by $P_i(1 \leq i \leq k)$.\\
\indent
Blue: $\{(i,i-k)~|~i \in\{2,4,\cdots,k\}\}$, $\{(i,i+k)~|~i \in\{1,3,\cdots,k-1\}\}$, $\{(i+(a-1)k,i+(a-2)k)~|~i\in\{2,4,\cdots,k\}\}$;\\
\indent
Purple: $\{(i,i-k)~|~i \in\{1,3,\cdots,k-1\}\}$, $\{(i,i+k)~|~i \in\{2,4,\cdots,k\}\}$, $\{(i+(a-1)k,i+(a-2)k),(n,n-k)~|~i \in\{3,5,\cdots,k-1\}\}$.

\vspace{-0.6em}
\begin{figure}[htbp]
\centering
\includegraphics[height=4.4cm, width=0.7\textwidth]{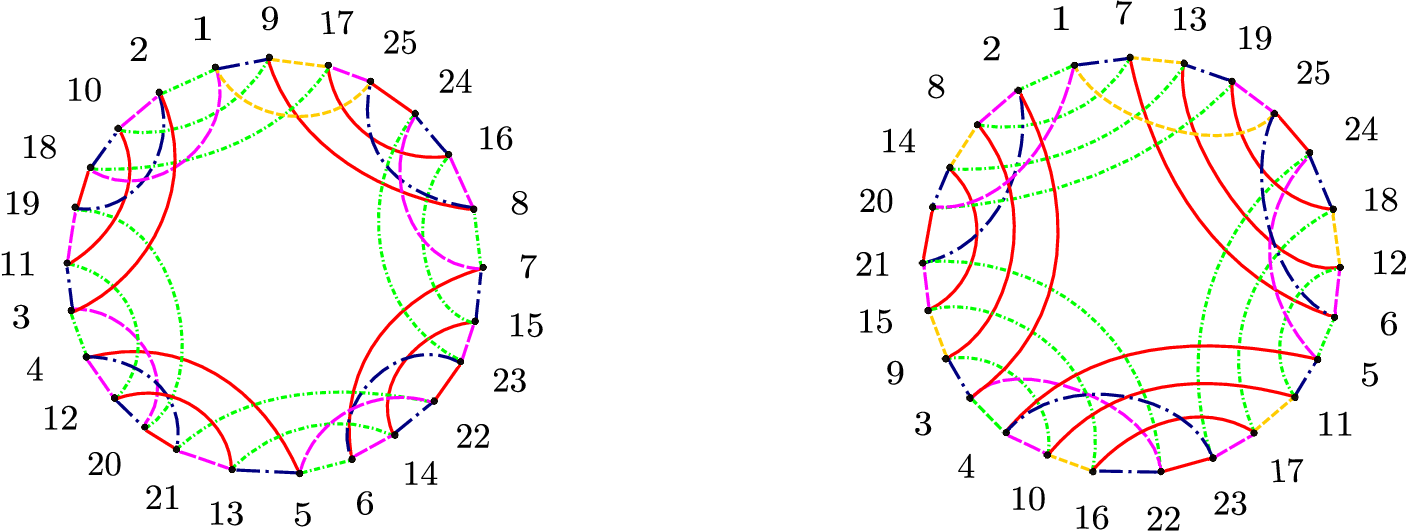}
\centerline{Fig.12  ~The matching book embedding of $C(25,8)$(left) and $C(25,6)$(right).}
\end{figure}
\vspace{-1.2em}

\noindent
\textbf{Step 3:} The coloring of the remaining edges induced by $P_i(1 \leq i \leq k)$.\\
\indent
As for all paths induced by $P_i(1 \leq i \leq k)$, two edges incident to the two endpoints have been colored with blue and purple respectively, then the remaining edges can be colored with blue and yellow alternately(see Fig.13(left)).

\vspace{-0.8em}
\begin{figure}[htbp]
\centering
\includegraphics[height=4.9cm, width=0.72\textwidth]{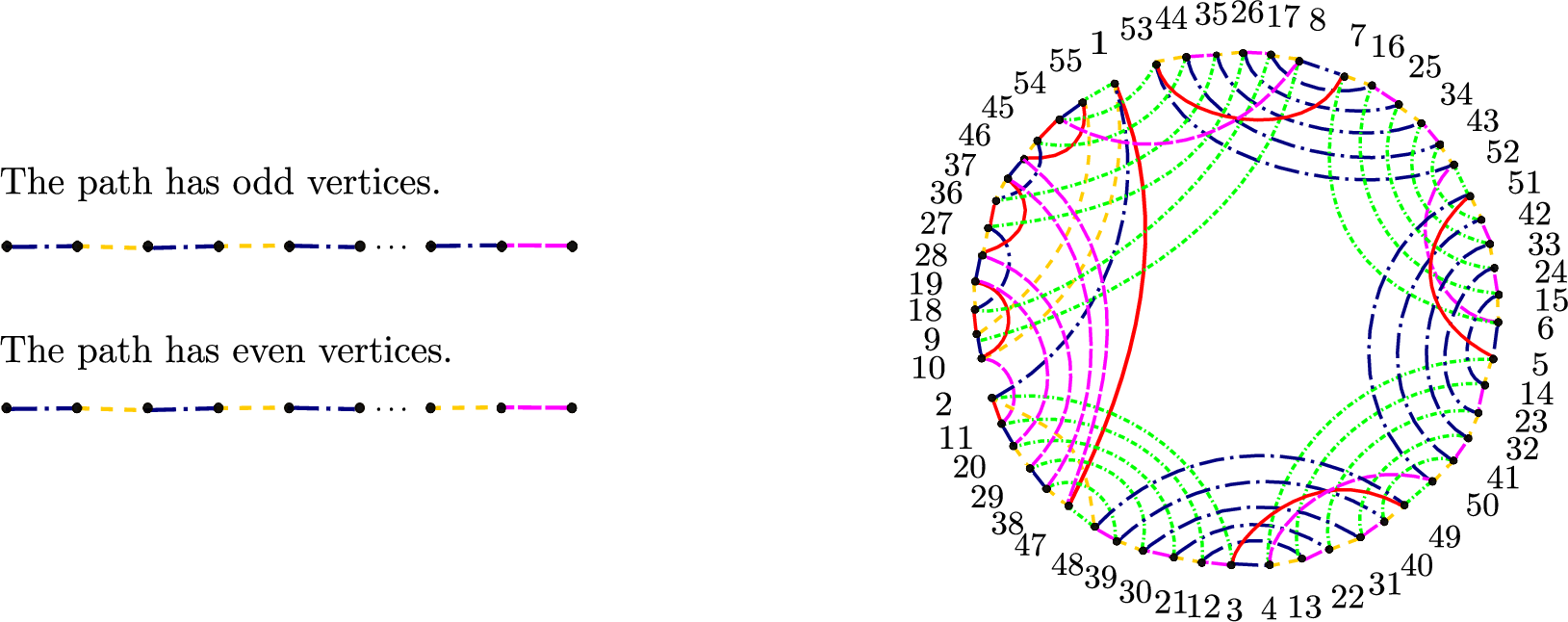}
\centerline{Fig.13  ~The coloring of path(left) and the matching book embedding of $C(55,9)$(right).}
\end{figure}
\vspace{-0.8em}

\vspace{0.3em}
\noindent\textbf{Subcase 2.2.} $k$ is odd.\\
\indent
Let $Q_1,Q_2,\cdots,Q_{k+1}$ be ordered vertex sets, where $Q_1$$=$$\{1\},Q_{i}$$=$$P_i(2 \leq i \leq k),Q_{k+1}$$=$$P_1\backslash\{1\}$. We denote the $j$-th element in $Q_i$ by $Q_{i_{j}}$.~Since $k$ is odd, $a(a \geq 4)$ and $|Q_i|(2 \leq i \leq k+1)$ are even.~Let $Q_c$ be an ordered vertex set which is a rearrangement of $Q_k\cup Q_{k+1}$, specifically, $Q_c$$=$$\{Q_{(k+1)_1},Q_{k_1},Q_{k_2},Q_{(k+1)_2},Q_{(k+1)_3},Q_{k_3},Q_{k_4},Q_{(k+1)_4},\cdots,Q_{(k+1)_{a-1}},Q_{k_{a-1}},Q_{k_a},Q_{(k+1)_a}\}.$

Put the vertices of $G$ counterclockwise along a circle in the order $Q_1Q_{c}^{-}Q_2Q_3^{-}\cdots Q_{k-2}^{-}Q_{k-1}$. The edges of $G$ can be colored well in the following way(see Fig.13(right) for $C(55,9)$).\\
\indent Green: $\{(i,i+1)~|~i \in Q_j, j\in\{2,4,\cdots,k-1\}\}$, $\{(1,n)\}$;\\
\indent Blue: $\{(i,i+1)~|~i \in Q_j, j\in \{1,3,\cdots,k-2\}\}$,~$\{(k,k+1),(n-1,n)\}$,~$\{(i+jk,i+(j+1)k)~|~i \in\{2,k,k+1\},j \in \{1,3,\cdots,a-3\} \}$;\\
\indent Purple: $\{(i,i+1)~|~i \in Q_{k+1}\backslash\{n\}\}$,
$\{(i,i-k)~|~i \in \{4,6,\cdots,k-1\}\}$, $\{(3+jk,3+(j+1)k)~|~j\in \{0,2,\cdots,a-2\}\}$, $\{(i+jk,i+(j+1)k)~|~i \in \{4,5,\cdots,k-1\},j\in \{1,3,\cdots,a-3\}\}$;\\
\indent Yellow: $\{(i,i+1)~|~i \in Q_k\backslash\{k,n-1\}\}$, $\{(1,1+k),(k,n),(2,2-k)\}$, $\{(2+jk,2+(j+1)k)~|~j \in\{2,4,\cdots,a-2\}\}$, $\{(3+jk,3+(j+1)k)~|~j\in \{1,3,\cdots,a-3\}\}$, $\{(i+jk,i+(j+1)k)~|~i \in \{4,5,\cdots,k-1\},j\in \{0,2,\cdots,a-2\}\}$;\\
\indent Red: $\{(2,2+k)\}$, $\{(i+jk,i+(j+1)k)~|~i \in\{k,k+1\},j\in\{0,2,\cdots,a-2\}\}$, $\{(i,i-k)~|~i \in \{1,3,\cdots,k-2\}\}$.\\
\indent
So $mbt(G)\leq5$. The result is established.
\hfill{$\square$}

\vspace{0.4em}
\noindent \textbf{Theorem 6.2.~}\emph{Let $G=C(n,k)$, where $n$ is odd. If~$r=k-1$, then $mbt(G)=\Delta(G)+1=5.$}\\
\noindent \textbf{Proof.}
By Lemma 2.5, it is sufficient to show that $mbt(G)\leq\Delta(G)+1$.

\vspace{0.2em}
\noindent\textbf{Case 1:} $k$ is even\\
\indent
Put all vertices of $G$ clockwise along a circle in the order $P_1P_2^{-}P_3P_4^{-}\cdots P_{k-1}P_k^{-}$.~The edges of $G$ can be colored well in the following two steps (see Fig.14(left) for $C(53,6)$):

\vspace{-0.4em}
\begin{figure}[htbp]
\centering
\includegraphics[height=5.2cm, width=0.75\textwidth]{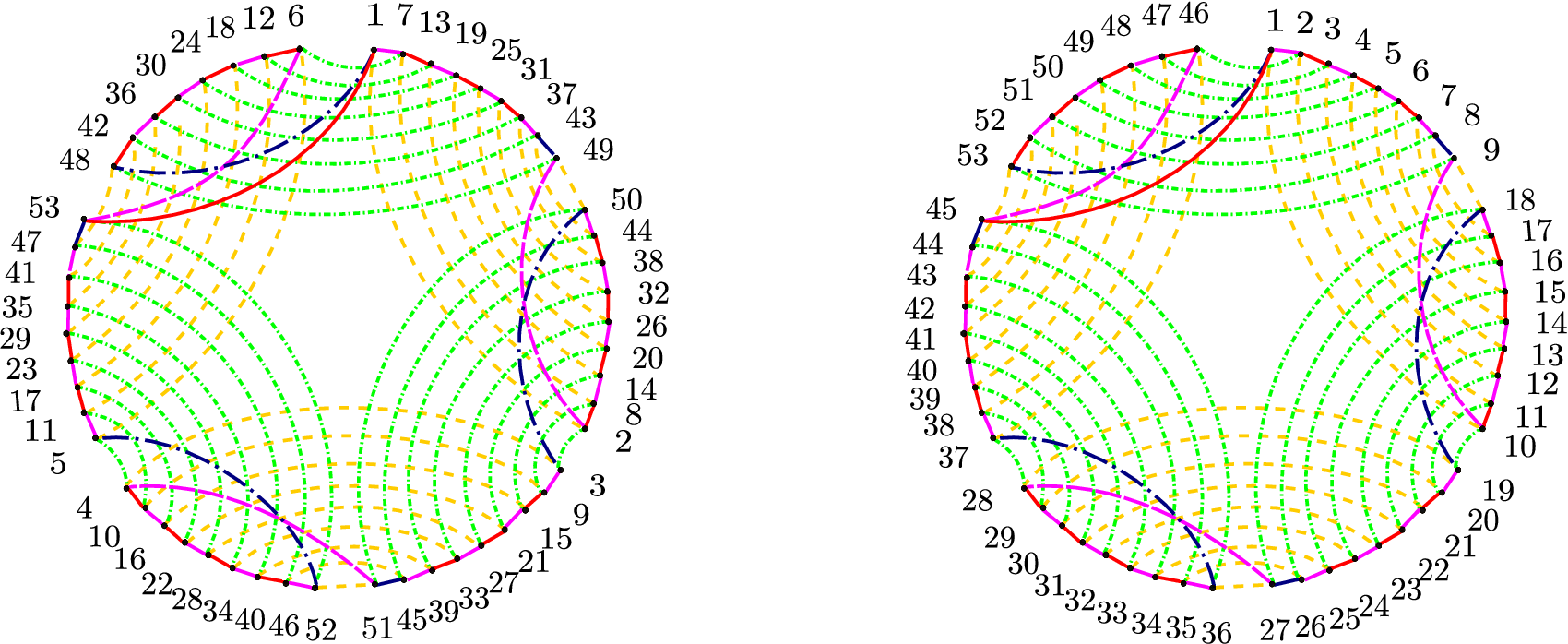}
\centerline{Fig.14  ~The matching book embedding of $C(53,6)$(left) and $C(53,9)$(right).}%
\end{figure}

\vspace{-0.8em}
\noindent
\textbf{Step 1:} The coloring of the $n$-cycle $C_0$.

Color the edge $(1,n)$ red, the edges of $\{(i,i+1)~|~i \in P_j\backslash\{n\},j \in \{1,3,\cdots,k-1\}\}$ yellow, the edges of $\{(i,i+1)~|~i \in P_j,j \in\{2,4,\cdots,k\}\}$ green.

\noindent
\textbf{Step 2:} The coloring of the $n$-cycle $C_1$.\\
\indent
Color the edges $(i,i-k)(i \in \{1,3,\cdots,k-1\})$ blue, the edges $(i,i-k)(i \in \{2,4,\cdots,k\})$ purple.~It is easy to use red, blue and purple to color the remaining edges induced by $P_i(1 \leq i \leq k)$.

\noindent\textbf{Case 2:} $k$ is odd

Since $n,k$ are odd, $a$ is odd. Let $Q_i=\{1+(i-1)k,2+(i-1)k,\cdots,ik\}(1\leq i\leq a),Q_{a+1}=\{ak+1,ak+2,\cdots,n\}$ be ordered vertex sets. Put all vertices of $G$ clockwise along a circle in the order $Q_1Q_2^{-}Q_3Q_4^{-}\cdots Q_{a}Q_{a+1}^{-}$.~The edges of $G$ can be colored well with $\Delta(G)+1$ colors in the following two steps (see Fig.14(right) for $C(53,9)$):

\noindent
\textbf{Step 1:} The coloring of the $n$-cycle $C_1$.

Color the edge $(1,1-k)$ red, the edges of $\{(i,i+k)~|~i \in Q_j\backslash\{1-k\},j \in \{1,3,\cdots,a\}\}$ yellow, the edges of $\{(i,i+k)~|~i \in Q_j,j \in\{2,4,\cdots,a+1\}\}$ green.



\noindent
\textbf{Step 2:} The coloring of the $n$-cycle $C_0$.

Color the edges $(1+(i-1)k,(i-1)k)(i \in \{1,3,\cdots,a\})$ blue, the edges $(1+(i-1)k,(i-1)k)(i \in \{2,4,\cdots,a+1\})$ purple.~It is easy to use red, blue and purple to color the remaining edges induced by $Q_i(1 \leq i \leq a+1)$.\\
\indent
So $mbt(G)\leq5$. The result is established.
\hfill{$\square$}

\vspace{0.4em}
If $2\leq r \leq k-2$, we will discuss the dispersability of $C(n,k)$ by the parity of $t$ in Section 2. Whether $t$ is odd or even, $|P_i|=|P_j|>|P_t|(1 \leq i < j \leq t-1)$.

\vspace{0.4em}
\noindent \textbf{Theorem 6.3.~}\emph{Let $G=C(n,k)$, where $n,t(t \geq 3)$ are both odd. If $2 \leq r \leq k-2$, then $mbt(G)=\Delta(G)+1=5.$}\\
\noindent \textbf{Proof.~}According to Lemma 2.5, it is sufficient to show that $mbt(G)\leq\Delta(G)+1$.\\
\indent Since $t$ is odd and $|P_1|=|P_2|=\cdots=|P_{t-1}|$, $|P_t|$ is odd.~Let $P_i^{\prime},P_i^{\prime\prime}(i \in \{1,t-1\})$ be ordered vertex sets, which keeps the order in $P_i$ and the set $P_{i}^{\prime}$ includes the first $|P_t|$ elements of $P_i$, $P_{i}^{\prime\prime}$$=$$P_i\backslash P_{i}^{\prime}$. Let the ordered vertex set $P_{t}^{\prime}$ be a rearrangement of  $P_{t-1}^{\prime}\cup P_{t}$, specifically, $P_{t}^{\prime}=\{t-1,t,t+k,t-1+k,t-1+2k,t+2k,t+3k,t-1+3k,\cdots,t-1+(|P_t|-3)k,t+(|P_t|-3)k,t+(|P_t|-2)k,t-1+(|P_t|-2)k,t-1+(|P_t|-1)k,t+(|P_t|-1)k\}.$

\vspace{0.4em}
\noindent\textbf{Case 1:} $t=3$

Put all vertices of $G$ counterclockwise along a circle in the order $P_1^{\prime}P_3^{\prime-}P_1^{\prime\prime}P_2^{\prime\prime-}$.~The edges of $G$ can be colored well with $\Delta(G)+1$ colors in the following two steps.

\vspace{0.3em}
\noindent\textbf{Subcase 1.1.} $|P_t| > \lfloor |P_1|/2 \rfloor$(see Fig.15(left) for $C(63,11)$).

\noindent
\textbf{Step 1:} The coloring of the $n$-cycle $C_0$.

Color the edge $(2,3)$ red,~the edge $(2+k,3+k)$ blue,~the edges of $\{(i,i-1)~|~i \in P_2\}$ purple,~the edges of $\{(i,i-1)~|~i \text{~is~the~last~$|P_1^{\prime}|-|P_1^{\prime\prime}|$~elements~in~$P_1^{\prime}$}\} \cup \{(i-1,i-2)~|~i \in P_1^{\prime\prime}\}$ yellow,~the edges of $\{(i,i-1),(i+1,i+2)~|~i \in P_1^{\prime\prime}\} \cup \{(i,i-1)~|~i \in \{3+2k,3+3k,\cdots,3+(|P_3|-|P_2^{\prime\prime}|-1)k\}\}$
green.

\noindent
\textbf{Step 2:} The coloring of the $n$-cycle $C_1$.\\
\indent Color the edge $(n-1-k,n-1)$ red, the edges of $\{(i,i-k)~|~i \in\{1,2,3,n\}\}$ blue. Color the path induced by $P_1^{\prime}$ with red and blue alternately, where the edge $(1,1+k)$ is red. Color the path induced by $P_1^{\prime\prime}$ with yellow and blue alternately, the path induced by $P_2^{\prime\prime}$ with yellow and red alternately, where the edge $(2-2k,2-k)$ is yellow. As for the paths induced by $P_2^{'},P_3$, color the edges of $\{(i,i+k)~|~i \in\{2,3\}\}$ green, other edges with red and blue alternately, where the edges of $\{(i+k,i+2k)~|~i \in\{2,3\}\}$ are red.

\vspace{-0.6em}
\begin{figure}[htbp]
\centering
\includegraphics[height=5.8cm, width=0.82\textwidth]{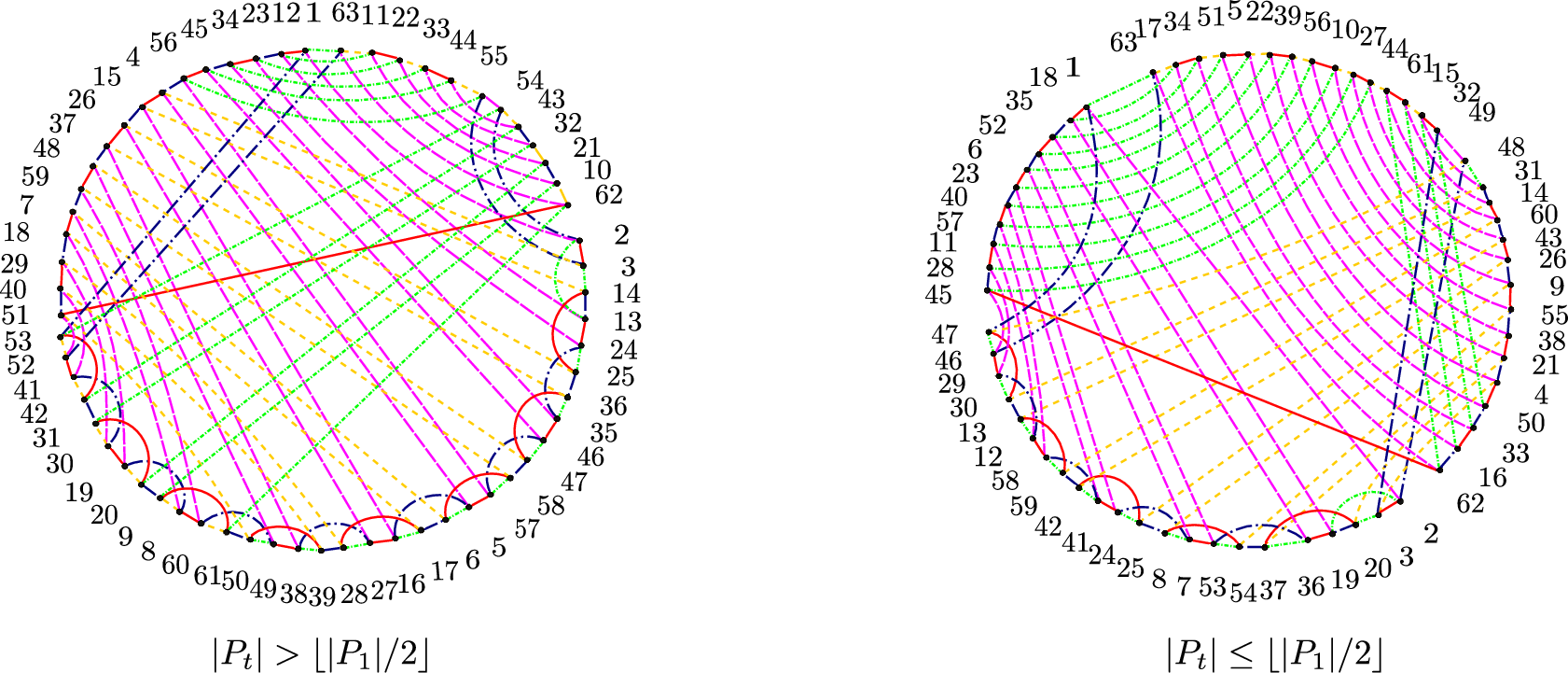}
\centerline{Fig.15  ~The matching book embedding of $C(63,11)$(left) and $C(63,17)$(right).}%
\end{figure}
\vspace{-1em}

\vspace{0.3em}
\noindent\textbf{Subcase 1.2.} $|P_t| \leq \lfloor |P_1|/2 \rfloor$(see Fig.15(right) for $C(63,17)$).\\
\noindent
\textbf{Step 1:} The coloring of the $n$-cycle $C_0$.\\
\indent Color the edge $(2,3)$ red,~the edge $(2+k,3+k)$ blue,~the edges of $\{(j,j-1)~|~j \in P_3\backslash\{3,3+k\}\} \cup \{(i,i-1)~|~i \text{~is~the first $|P_1|-|P_3|$ elements in~$P_1$}\}$ green, the edges of $\{(i,i-1)~|~i \in P_2\}$ purple, the edges of $\{(i,i-1)~|~i\text{ is the last } \text{$|P_3|$ elements in $P_1$}\}$ yellow.\\
\noindent\textbf{Step 2:} The coloring of the $n$-cycle $C_1$.\\
\indent Color the edge $(n-1-k,n-1)$ red, the edges of $\{(i,i-k)~|~i \in\{1,2,3,n\}\}$ blue, the edge $(2-k,2-2k)$ green. Color the path induced by $P_1^{\prime}$ with red and blue alternately, where the edge $(1,1+k)$ is red. Color the path induced by $P_1^{\prime\prime}\backslash\{2-k\}$ with red and blue alternately, where the edge $(n-1,k-1)$ is blue. Color the path induced by $P_2^{\prime\prime}$ with yellow and red alternately. As for the paths induced by $P_2^{'},P_3$, color the edges of $\{(i,i+k)~|~i \in\{2,3\}\}$ green, other edges with red and blue alternately, where the edges of $\{(i+k,i+2k)~|~i \in\{2,3\}\}$ are red.

\vspace{0.4em}
\noindent\textbf{Case 2: $t>3$}\\
\indent
Put all vertices of $G$ counterclockwise along a circle in the order $P_1^{-}P_2P_3^{-}P_4\cdots P_{t-4}^{-}P_{t-3}P_{t-1}^{\prime\prime}\\P_{t-2}^{-}P_{t}^{\prime}.$ All edges can be colored well with $\Delta(G)+1$ colors in the following four steps(see Fig.16 for $C(65,14)$(left) and $C(67,47)$(right)):%

\noindent
\textbf{Step 1: }The coloring of the $n$-cycle $C_0$.

Color the edge $(t-1,t)$ red, the edges of $\{(j,j+1)~|~j \in P_i, i\in\{1,3,\cdots,t-2\}\}$ yellow, the edges of $\{(j,j+1)~|~j \in P_i, i\in \{2,4,\cdots,t-3\} \cup \{t\} \}$ green, the edges of $\{(n-k,1-k)\}\cup\{(j,j+1)~|~j \in P_{t-1}^{\prime\prime}\}$ purple, the edges of $\{(i,i-1)~|~i\in P_t\backslash\{t,1-k\}\}$ blue.

\noindent
\textbf{Step 2: }The coloring of the edges connected the path induced by $P_i,P_j$$(i\equiv j+1(\text{mod} ~t))$.

Red:~$\{(i,i-k)~|~i \in \{2,4,\cdots,t-3\}\cup\{t-2\}\}$;

Blue:~$\{(i,i-k)~|~i \in \{1,3,\cdots,t-4\}\cup\{t-1,t\}\}$.

\noindent
\textbf{Step 3: }The coloring of the paths induced by $P_{t-1}^{\prime},P_t$.

Purple: $\{(i-1,i-1-k),(i,i-k)~|~i \text{~is the element in the even position of}~P_t\}$;

Red:~$\{(i,i+k),(i-1,i-1+k)~|~i \text{~is the element in the even position of}~ P_t\}$.

\vspace{-0.6em}
\begin{figure}[htbp]
\centering
\includegraphics[height=5.6cm, width=0.8\textwidth]{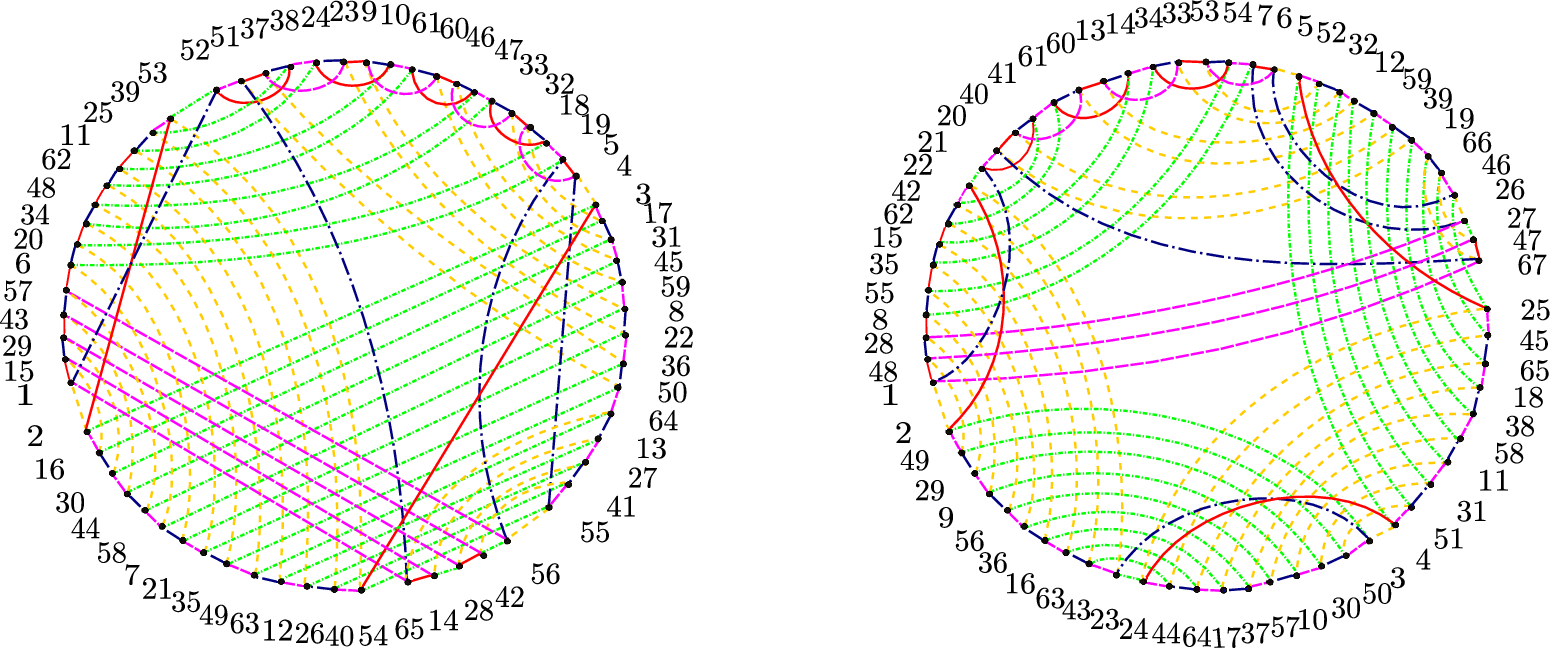}
\centerline{Fig.16  ~The matching book embedding of $C(65,14)$(left) and $C(67,47)$(right).}%
\end{figure}
\vspace{-1em}

\noindent
\textbf{Step 4:} The coloring of the paths induced by $P_i(i=1,2,\cdots,t-2)$, $P_{t-1}^{\prime\prime}\cup\{n-k\}$.

Color the path induced by $P_i(2\leq i \leq t-2)$ with purple and blue alternately.
~Color the edge $(1+(|P_1|-2)k,1+(|P_1|-1)k)$ purple, the remaining edges induced by $P_1$ with red and blue alternately, where the edge $(1,1+k)$ is red.
Color the edge $(n-k,n)$ blue, the path induced by $P_{t-1}^{\prime\prime}$ with green and red alternately.\\
\indent All edges of $G$ can be matching book embedded in five pages. Hence $mbt(G)=5$.
\hfill{$\square$}

\vspace{0.4em}
\noindent \textbf{Theorem 6.4.~}\emph{Let $G=C(n,k)$, where $n$ is odd and $t(t \geq 2)$ is even. If $2 \leq r \leq k-2$, then $mbt(G)=\Delta(G)+1=5.$}

\noindent \textbf{Proof.}
By Lemma 2.5, it is sufficient to show that $mbt(G)\leq \Delta(G)+1$.

\vspace{0.4em}
\noindent\textbf{Case 1: $t=2$}\\
\indent
Let $P_j^{\prime}$, $P_j^{\prime\prime}(j \in\{1,2\})$ be ordered subsets of $P_j$, which keeps the order of $P_j$, where $P_j^{\prime}$=$\{~i~|~i \in P_j, i \text{~is~the~element~before~the~element~$n+1-j$~in~$P_j$}\}, P_j^{\prime\prime}$=$P_j\backslash P_j^{\prime}$.
Let the ordered vertex set $L$ be a rearrangement of $P_1^{\prime\prime}\cup P_2^{\prime\prime}$, specifically, $L$$=$$\{n,n-1,k-1,k,2k,2k-1,3k-1,3k,\cdots,(|P_1^{\prime\prime}|-3)k,(|P_1^{\prime\prime}|-3)k-1,(|P_1^{\prime\prime}|-2)k-1,(|P_1^{\prime\prime}|-2)k,(|P_1^{\prime\prime}|-1)k,(|P_1^{\prime\prime}|-1)k-1\}$.

Put the vertices of $G$ counterclockwise along a circle in the order $P_1^{\prime}L^{-}P_2^{\prime-}$. The edges of $G$ can be colored well with $\Delta(G)+1$ colors in the following three steps(see Fig.17(left) for $C(47,13)$):

\noindent
\textbf{Step 1: }The coloring of the $n$-cycle $C_0$.

Color the edges of $\{(i-1,i)~|~i \in P_2^{\prime}\}\cup\{(i,i+1)~|~i \in P_2^{\prime\prime}\backslash\{n-1\}\}$ purple, the edges of  $\{(1,n)\}\cup\{(i,i-1)~|~i \in P_2^{\prime\prime}\}$ red, the edges of  $\{(n-1,n)\}\cup\{(i,i+1)~|~i \in P_2^{\prime}\}$ green, the edges of $\{(i,i+1)~|~i \in P_1^{\prime\prime}\backslash\{n\}\}$ yellow.




\noindent
\textbf{Step 2:} The coloring of the paths induced by $P_1^{\prime\prime}\cup\{2,n-k\},P_2^{\prime\prime}\cup\{1,n-1-k\}$.\\
\indent Color the edge $(n-k,n)$ purple, the edge $(n-1,n-1-k)$ yellow, the edges of  $\{(i,i+k),(i-1,i-1+k)~|~i\text{~is~the~element~in~the~even~position~of~$P_1^{\prime\prime}$}\}$ green, and the edges of $\{(i,i+k),(i-1,i-1+k)~|~i \text{~is~the~element~in~the~odd~position~of~$P_1^{\prime\prime}$}\}$ blue.\\
\noindent
\textbf{Step 3:} The coloring of the paths induced by $P_1^{\prime}$ and $P_2^{\prime}$.

The edges of the path induced by $P_2^{\prime}$ can be colored with yellow and red alternately, where the edge $(n-1-2k,n-1-k)$ is red. Because $n$ is odd and $t=2$, $|P_1^{\prime\prime}|$ is odd and equals to $|P_2^{\prime\prime}|$. There are odd elements in $\{i~|~i+2 \in P_1^{\prime\prime}\}$, which is a subset of $P_1^{\prime}$ and induces a path with even edges. Color these edges with purple and blue alternately, where the edge $(k-2,n-2)$ is purple. Color the edge $(n-k-2,n-2)$ blue. The remaining edges induced by $P_1^{\prime}$ can be colored with green and red alternately, where the edge $(1,1+k)$ is green.

\vspace{-0.6em}
\begin{figure}[htbp]
\centering
\includegraphics[height=5.6cm, width=0.82\textwidth]{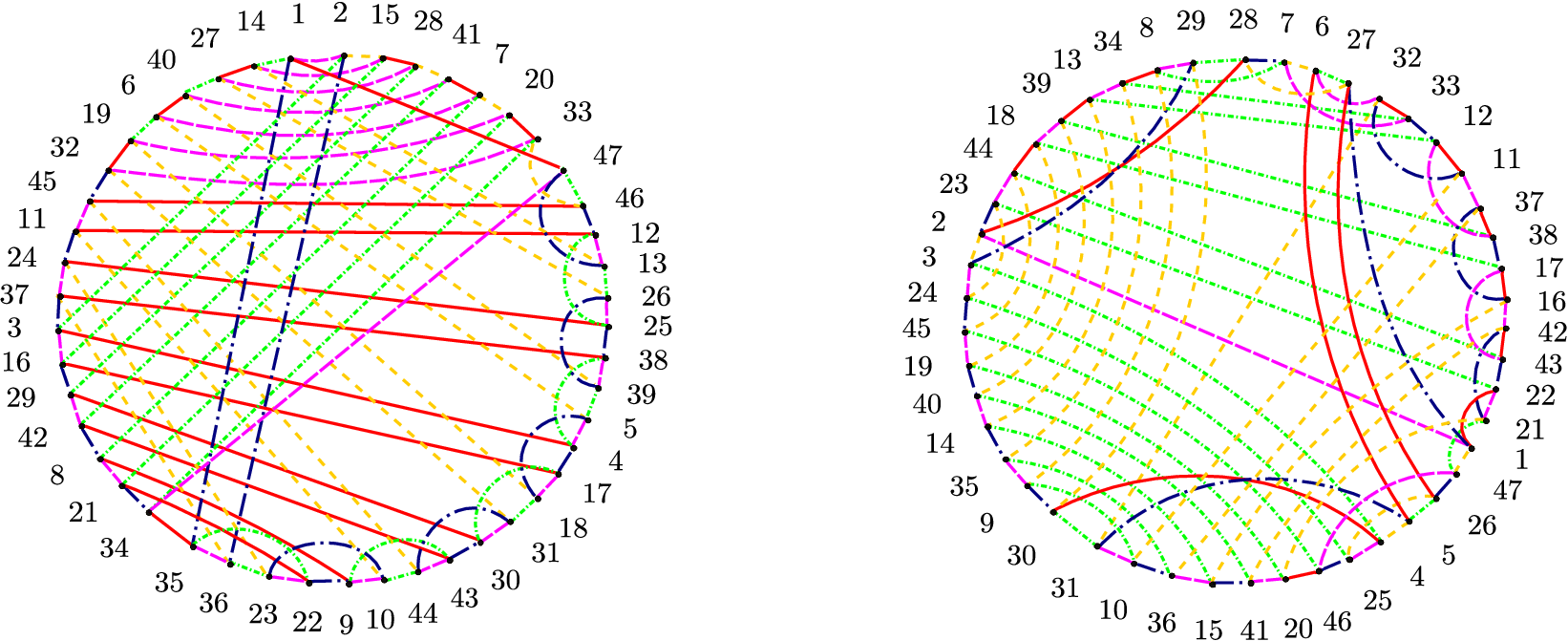}
\centerline{Fig.17  ~The matching book embedding of $C(47,13)$(left) and $C(47,21)$(right).}%
\end{figure}
\vspace{-1em}

\vspace{0.3em}
\noindent\textbf{Case 2: $t\geq4$}\\
\indent
Let $P_j^{\prime},P_j^{\prime\prime}(j \in \{1,t-1\})$ be ordered vertex subsets of $P_j$,~which keeps the order in $P_j$,~where $P_1^{\prime}$=$\{~i~|~i \text{~is~the~element~before~the~element~$t+1$~in~$P_{1}$}\}$, 
$P_{t-1}^{\prime}=\{~i~|~i \text{~is~the~element~before~the}\\ \text{element~$n$~in~$P_{t-1}$}\}$ and $P_{j}^{\prime\prime}$=$P_j\backslash P_{j}^{\prime}$. 
Let the ordered set $L$ be a rearrangement of $P_1^{\prime}\cup P_{t-1}^{\prime\prime}$, where $L=\{n,1,k,1+k,1+2k,2k,3k,1+3k,1+4k,4k,\cdots,(|P_1^{\prime}|-2)k,1+(|P_1^{\prime}|-2)k,1+(|P_1^{\prime}|-1)k,(|P_1^{\prime}|-1)k\}$.
Put the vertices of $G$ counterclockwise along a circle in the order $P_2^{-}P_3P_4^{-}P_5\cdots $ $P_{t-2}^{-}P_{t-1}^{\prime}LP_{t}^{-}P_1^{\prime\prime}$.~The edges of $G$ can be colored well with $\Delta(G)+1$ colors in the following four steps(see Fig.17(right) for $C(47,21)$):

\noindent
\textbf{Step 1: }The coloring of the $n$-cycle $C_0$.

Yellow: $\{(j,j+1)~|~j \in P_2 \cup P_4 \cdots \cup P_{t-2}  \cup P_{t}\backslash\{n-1\}\}$, $\{(1,n)\}$;

Green: $\{(j,j+1)~|~j \in P_{1}\cup P_3 \cup \cdots \cup P_{t-3}\backslash\{1\} \}$;

Purple: $\{(1,2),(n-1,n),(k,k+1)\}$;

Red: $\{(j+1,j)~|~j \in P_{t-1}\backslash\{n,k\}\}$.

\noindent
\textbf{Step 2:} The coloring of the paths induced by $P_1^{\prime}\cup\{t+1\},P_{t-1}^{\prime\prime}\cup\{t,n-k\}$.

Green: $\{(k,n)\}$;

Blue: $\{(j-1,j+k-1),(j,j+k)~|~\text{$j$~is~the~element~in~the~even~position~of~$P_{1}^{\prime}$}\}$, $\{(n-k,n)\}$;

Purple: $\{(j-1,j+k-1),(j,j+k)~|~\text{$j$~is~the~element~in~the~even~position~of~$P_{1}^{\prime}\backslash\{1\}$}\}$;

Red: $\{(1,1+k)\}$.

\noindent
\textbf{Step 3: }The coloring of the edges connected the paths induced by $P_i,P_j$$(i\equiv j+1(\text{mod} ~t),i \neq t)$.

Blue: $\{(i,i-k)~|~ i \in \{1,3,\cdots,t-1\}\}$;

Red: $\{(i,i-k)~|~ i \in \{2,4,\cdots,t-2\}\}$.

\noindent
\textbf{Step 4:} The coloring of the edges of the paths induced by $P_1^{\prime\prime}$, $P_i(2 \leq i \leq t-2)$, $P_{t-1}^{\prime}$, $P_t$.

The edges of the path induced by $P_1^{\prime\prime}$ can be colored with blue and purple alternately, where the edge $(t+1,t+1+k)$ is blue. The edges of the path induced by $P_{t-1}^{\prime}, P_t$ can be colored with green and purple alternately, where the edge $(t,t+k)$ is green. As for $3\leq i \leq t-3$, all edges of the path induced by $P_i$ can be colored with blue and purple alternately. Color the edge $(2,2+k)$ blue, and use red and purple to color the remaining edges of $P_{2}$. Color the edge $(n-1-k,n-1)$ blue, the edge $(k-1,n-1)$ red, and use purple and blue to color the remaining edges induced by $P_{t-2}$.

Hence $mbt(G)\leq 5$. The result is established.
\hfill{$\square$}

\section{Conclusion}
The first conclusion to be drawn is that the classification of the circulant graphs $C(\mathbb{Z}_{n},\{ k_1, k_2\})$ is obtained mainly by the Diophantine equation technique. Secondly,
 the dispersability of $G=C(\mathbb{Z}_{n},\{ 1, k\})$ for any $2\leq k \leq n-2$ is determined. Specifically, if $G$ is a bipartite graph, then $G$ is dispersable; if $G$ contains an odd cycle, then  $G$ is nearly dispersable.
 A different method can be used to get the dispersability of the Cartesian graph bundle over cycles, and the resulting paper will be published elsewhere.


\section*{Acknowledgment}

\indent

We would like to thank an anonymous referee for very helpful comments and suggestions. This work was partially funded by Science and Technology Project of Hebei Education Department, China (No.~ZD2020130) and the Natural Science Foundation of Hebei Province, China (No. A2021202013).

\end{document}